\documentclass[a4paper,
               10pt,
               pdftex,
               normalheadings,
               headsepline,
               footsepline,
               headinclude,
               footinclude,
               DIV14,
               abstracton]
{scrartcl}

\usepackage[left=2.50cm, right=2.50cm, top=2.00cm, bottom=3.50cm]{geometry}

\usepackage
{
    graphicx,
    amssymb,
    amsmath,
    amsthm,
    xcolor,
    dsfont,
    algpseudocode,
    authblk,
}
\usepackage[ngerman, english]{babel}

\usepackage{enumerate}
\usepackage{algorithm} 
\usepackage{algpseudocode}
\usepackage{tabularx}
\usepackage{subcaption} 
\usepackage{booktabs} 
\usepackage{algpseudocode}
\usepackage{tikz}
\usetikzlibrary{decorations.pathreplacing}
\usepackage{caption}
\usepackage{subcaption}
\usepackage{bibentry}

\usepackage[colorlinks,
            pdffitwindow=false,
            plainpages=false,
            pdfpagelabels=true,
            pdfpagemode=UseOutlines,
            pdfpagelayout=SinglePage,
            bookmarks=false,
            colorlinks=true,
            hyperfootnotes=false,
            linkcolor=black,
            citecolor=black]
{hyperref}

\usepackage{mathtools}
\usepackage{enumerate}
\usepackage{thmtools}
\usepackage{cite}
\usepackage{float}

\usepackage{sectsty}
\sectionfont{\fontsize{10}{10}\selectfont}
\subsectionfont{\fontsize{10}{10}\selectfont}

\newcommand{\R}{\mathbb{R}}

\renewcommand{\H}{\mathcal{H}}
\newcommand{\W}{\mathcal{W}}
\newcommand{\E}{\mathcal{E}}

\newcommand{\logeq}{\mathrel{\raisebox{.66pt}{:}}\Leftrightarrow}

\DeclareMathOperator*{\argmin}{arg\,min}

\DeclareMathOperator*{\proj}{proj}

\DeclareMathOperator{\co}{conv}

\makeatletter
\newcommand{\leqnomode}{\tagsleft@true\let\veqno\@@leqno}
\newcommand{\reqnomode}{\tagsleft@false\let\veqno\@@eqno}
\makeatother

\newtheorem{theorem}{Theorem}[section]
\newtheorem{myprop}[theorem]{Proposition}
\newtheorem{mydef}[theorem]{Definition}
\newtheorem{mycorollary}[theorem]{Corollary}
\newtheorem{mylemma}[theorem]{Lemma}
\newtheorem{myremark}[theorem]{Remark}
\newtheorem{myassump}[theorem]{Assumption}
\newtheorem{myexample}[theorem]{Example}

\newcommand{\ignore}[1]{}
\newcommand{\nobibentry}[1]{{\let\nocite\ignore\bibentry{#1}}}

\let\OLDthebibliography\thebibliography
\renewcommand\thebibliography[1]{
	\OLDthebibliography{#1}
	\setlength{\parskip}{0pt}
	\setlength{\itemsep}{0pt plus 0.3ex}
}

\begin{document}

\title{Fast Convergence of Inertial Multiobjective Gradient-like Systems with Asymptotic Vanishing Damping}
\author[1]{Konstantin Sonntag}
\author[2]{Sebastian Peitz}
\affil[1]{\normalsize Department of Mathematics, Paderborn University, Germany}
\affil[2]{\normalsize Department of Computer Science, Paderborn University, Germany}

\date{}

\maketitle


\begin{abstract}
We present a new gradient-like dynamical system related to unconstrained convex smooth multiobjective optimization which involves inertial effects and asymptotic vanishing damping. To the best of our knowledge, this system is the first inertial gradient-like system for multiobjective optimization problems including asymptotic vanishing damping, expanding the ideas previously laid out in [\textsc{H. Attouch and G. Garrigos}, \emph{Multiobjective optimization: an inertial dynamical approach to pareto optima}, arXiv preprint arXiv:1506.02823, (2015)]. We prove existence of solutions to this system in finite dimensions and further prove that its bounded solutions converge weakly to weakly Pareto optimal points. In addition, we obtain a convergence rate of order $\mathcal{O}(t^{-2})$ for the function values measured with a merit function. This approach presents a good basis for the development of fast gradient methods for multiobjective optimization.
\end{abstract}

\section{Introduction}
\label{sec:introduction}
In this paper $\H$ is a real Hilbert space with scalar product $\langle \cdot, \cdot \rangle$ and norm $\lVert \cdot \rVert$. We are interested in a gradient-dynamic approach to the Pareto optima of the multiobjective optimization problem
\leqnomode
\begin{align}
\label{eq:MOP}
    \min_{x \in \H} F(x),
\tag{MOP}
\end{align}
with $F: \H \to \R^m, F(x) = \left(f_1(x), \dots, f_m(x)\right)$, where $f_i:\H \to \R$ are convex and continuously differentiable functions for $i = 1,\dots, m$. We define the following \emph{multiobjective inertial gradient-like dynamical system with asymptotic vanishing damping}
\begin{align}
\label{eq:MAVD}
\tag{MAVD}
        \frac{\alpha}{t}\dot{x}(t) + \proj_{C(x(t)) + \ddot{x}(t)}{(0)} = 0,
\end{align}
where $C(x) \coloneqq \co\left( \lbrace \nabla f_i(x) : i = 1,\dots, m\rbrace \right)$ is defined as the convex hull of the gradients. For a closed convex set $K \subset \H$ and a vector $x \in \H$, the projection of $x$ on the set $K$ is denoted by $\proj_K(x) \coloneqq \argmin_{y \in K} \lVert y - x \rVert^2$. Our interest in the system \eqref{eq:MAVD} is motivated by the active research in dynamical systems for fast minimization and their relationship with numerical optimization methods.

In the singleobjective case $m = 1$, \eqref{eq:MAVD} reduces to the \emph{inertial gradient system with asymptotic vanishing damping}
\begin{align*}
    \label{eq:AVD}
        \tag{AVD}
        \ddot{x}(t) + \frac{\alpha}{t}\dot{x}(t) + \nabla f(x(t)) = 0,
\end{align*}
which is introduced in \cite{Su2014} in connection with \emph{Nesterov's accelerated gradient method} (see \cite{Nesterov1983}) and analyzed further in \cite{Cabot2009_long, Cabot2009, Attouch2018}. For $\alpha > 0$ every solution $x$ of \eqref{eq:AVD} satisfies $\lim_{t \to + \infty} f(x(t)) = \min_{x \in \H} f(x)$. For $\alpha \ge 3$ it holds that $f(x(t)) - \min_{x\in \H} f(x) = \mathcal{O}(t^{-2})$ (see \cite{Su2014}). For $\alpha > 3$ the trajectories experience an improved converge rate of order $f(x(t)) - \min_{x\in \H} f(x) = o(t^{-2})$ and every solution $x$ converges weakly to a minimizer of $f$ given that the set of minimizers is nonempty (see \cite{Attouch2018, May2017}). Here, for a real valued function $g:[t_0 ,+\infty) \to \R_{\ge 0}$, with $t_0 > 0$, we write $g(t) = \mathcal{O}(t^{-2})$ if there exists $C > 0$ such that $t^2g(t) \le C$ for all $t \ge t_0$ and we write $g(t) = o(t^{-2})$ if $\lim_{t \to +\infty} t^2g(t) =0$. It is an open question if similar results can be obtained for multiobjective optimization problems (see \cite{Attouch2015}).

While there exists exhaustive literature on gradient-systems connected with singleobjective optimization problems, similar systems for multiobjective optimization problems are rarely adressed in the literature. There are only few results in this area which we present in the following. 

When we neglect the inertial effects introduced by $\ddot{x}(t)$ and drop the damping coefficient $\frac{\alpha}{t}$ in \eqref{eq:MAVD}, we return to the \emph{multiobjective gradient system}
\begin{align}
\label{eq:MOG}
        \tag{MOG}
        \dot{x}(t) + \proj_{C(x(t))}(0) = 0.
\end{align}
The relation of the system \eqref{eq:MOG} to multiobjective optimization problems is discussed in \cite{Miglierina2004, Attouch2014}. In \cite{Attouch2015_3} the system \eqref{eq:MOG} is extended to the setting of nonsmooth multiobjective optimization. A remarkable property of the system \eqref{eq:MOG} is that the function values of each objective decrease along any solution $x$ of \eqref{eq:MOG}, i.e., $\frac{d}{dt}f_i(x(t)) \le 0$ for all $t \ge t_0$ and for all $i = 1,\dots, m$. Further, bounded solutions of \eqref{eq:MOG} converge weakly to weakly Pareto optimal solutions (see \cite{Attouch2014}).

A first study on inertial gradient-like dynamical systems for multiobjective optimization in Hilbert spaces is proposed in \cite{Attouch2015}. The authors of \cite{Attouch2015} combine the system \eqref{eq:MOG} with the so-called \emph{heavy ball with friction dynamic}. For a scalar optimization problem $\min_{x\in \H} f(x)$ with a smooth objective function $f:\H \to \R$ the heavy ball with friction dynamical system reads as
\begin{align}
\label{eq:HBF}
        \tag{HBF}
        \ddot{x}(t) + \alpha \dot{x}(t) + \nabla f(x(t)) = 0,
\end{align}
for $\alpha > 0$. The system \eqref{eq:HBF} and its connection with optimization problems is well-studied (see \cite{Polyak1964, Alvarez2000, Attouch2000, Goudou2009}). It can be shown that for a convex smooth function $f$ with a nonempty set of minimizers, the solutions $x$ of \eqref{eq:HBF} converge weakly to minimizers of $f$ and $\lim_{t \to + \infty} f(x(t)) = \min_{x \in \H} f(x)$ (see \cite{Attouch2000}).

In \cite{Attouch2015} the systems \eqref{eq:MOG} and \eqref{eq:HBF} are combined to the \emph{inertial multiobjective gradient system}
\begin{align}
\label{eq:IMOG}
        \tag{IMOG}
        \ddot{x}(t) + \alpha \dot{x}(t) + \proj_{C(x(t))}(0) = 0,
\end{align}
for $\alpha > 0$ (see also \cite{Garrigos2015}). Any bounded solution $x$ of \eqref{eq:IMOG} converges weakly to a weak Pareto optimum of \eqref{eq:MOP} under the condition $\alpha^2 > L$, where $L > 0$ is a common Lipschitz constant of the gradients $\nabla f_i$. This last condition is the reason why it is not straight-forward to introduce asymptotic vanishing damping in \eqref{eq:IMOG}. If we include a damping coefficient $\frac{\alpha}{t}$ in the system \eqref{eq:IMOG} to adapt the proof, we would need $\left(\frac{\alpha}{t}\right)^2 > L$ for all $t \ge t_0$ which cannot hold. So either one has to find a different proof for the convergence of the trajectories of \eqref{eq:IMOG} or one has to define another generalization of the system \eqref{eq:HBF} to the multiobjective optimization setting.

In a first step to define a gradient-like system with asymptotic vanishing damping, in \cite{Sonntag2022} we defined the system
\begin{align}
\label{eq:IMOG'}
        \tag{IMOG'}
        \alpha \dot{x}(t) + \proj_{C(x(t)) + \ddot{x}(t)}(0) = 0,
\end{align}
with $\alpha > 0$, which also simplifies to the system \eqref{eq:HBF} in the singleobjective case. In \cite{Sonntag2022} it is shown that for smooth convex objective functions $f_i$ each bounded solution of \eqref{eq:IMOG'} converges weakly to a weak Pareto optimum of \eqref{eq:MOP}. When we introduce asymptotic vanishing damping in \eqref{eq:IMOG'}, we recover the system \eqref{eq:MAVD} which is analyzed in this paper.

In singleobjective optimization, asymptotic vanishing damping is of special interest from an optimization point of view since it guarantees fast convergence rates of the function values, namely of order $\mathcal{O}(t^{-2})$. To prove these convergence rates, one uses Lyapunov type energy functions which usually involve terms of the form $t^2(f(x(t)) - f(x^*))$, where $x^*$ is a minimizer of $f$. The choice of the minimizer is not crucial for the derivation since all minimizers have the same function value for a convex problem. In multiobjective optimization this is not possible anymore. Since there is no total order on $\R^m$, there is not a single solution to the multiobjective optimization problem \eqref{eq:MOP} but a set of solutions, the Pareto set. The values of the different objective functions vary along the Pareto front. If we choose a Pareto optimal solution $x^*$ to \eqref{eq:MOP} we run into the problem that the terms $f_i(x(t)) - f_i(x^*)$ do not necessarily remain positive along the solution trajectories of \eqref{eq:MAVD}. We need a different concept to define suitable energy functions to ensure monotonous decay of the energy along trajectories. A fruitful approach to analyze the convergence rate of multiobjective optimization methods is conducted by the use of so-called merit functions. In \cite{Tanabe2022_3}, the function $u_0(x) = \sup_{z \in \H} \min_{i = 1,\dots, m} f_i(x) - f_i(z)$ gets introduced as a merit function for nonlinear multiobjective optimization problems. The investigation of merit functions is part of the active research on accelerated first order methods for multiobjective optimization (see \cite{ElMoudden2021, Tanabe2022_2, Tanabe2022_4}). The function $u_0$ is nonnegative, attains the value zero only for weakly Pareto optimal solutions and is lower semicontinuous. It is therefore suitable as a measure of convergence speed for optimization methods in multiobjective optimization. In addition, in the singleobjective setting $u_0(x)$ simplifies to $f(x) - \inf_{z \in \H} f(z)$. Using this idea, we are able to define energy functions in the spirit of \cite{Su2014} to prove convergence rates of order $u_0(x(t)) = \mathcal{O}(t^{-2})$ for solutions $x$ to \eqref{eq:MAVD}.

This paper is organized as follows. In Section \ref{sec:moo}, we present the background on multiobjective optimization and merit functions which we use in our convergence analysis. We prove the existence of global solutions to the system \eqref{eq:MAVD} in finite dimensions in Section \ref{sec:global_existence}. Section \ref{sec:asym_behaviour_traj} contains our main results on the properties of the trajectories $x$ of \eqref{eq:MAVD}. In Theorem \ref{thm:conv_W}, we prove $\lim_{t \to + \infty} u_0(x(t)) + \frac{1}{2}\lVert \dot{x}(t) \rVert^2 = 0$ for all $\alpha > 0$. Convergence of order $u_0(x(t)) = \mathcal{O}(t^{-2})$ for $\alpha \ge 3$ is proven in Theorem \ref{thm:u_0_x(t)_=_O(t^-2)}. Theorem \ref{thm:weak_conv} states weak convergence of the bounded solutions of \eqref{eq:MAVD} to weakly Pareto optimal solutions using Opial's lemma, given $\alpha > 3$. In Section \ref{sec:numerical_experiments}, we verify the bounds for the convergence speed of $u_0(x(t))$ on two numerical examples. We present the conclusion and outlook on future research in Section \ref{sec:conclusion}.

\reqnomode
\section{Multiobjective optimization}
\label{sec:moo}
\subsection{Pareto optimal and Pareto critical points}
\label{subsec:Parerto_optimal_and_Pareto_ciritical_points}
The goal of multiobjective optimization is to optimize several functions simultaneously. In general, it is not possible to find a point minimizing all objective functions at once. Therefore, we have to adjust the definition of optimality in this setting. This can be done via the concept of Pareto optimality (see \cite{Miettinen1998}), which is defined as follows.
\begin{mydef}
Consider the optimization problem \eqref{eq:MOP}.
\begin{enumerate}[(i)]
    \item A point $x^* \in \H$ is \emph{Pareto optimal} if there does not exist another point $x \in \H$ such that $f_i(x) \le f_i(x^*)$ for all $i = 1,\dots,m$ and $f_j(x) < f_j(x^*)$ for at least one index $j$. The set of all Pareto optimal points is the \emph{Pareto set}, which we denote by $P$. The set $F(P) \subset \R^m$ in the image space is called the \emph{Pareto front}.
    \item A point $x^* \in \H$ is \emph{weakly Pareto optimal} if there does not exist another vector $x \in \H$ such that $f_i(x) < f_i(x^*)$ for all $i = 1,\dots, m$. The set of all weakly Pareto optimal points is the \emph{weak Pareto set}, which we denote by $P_w$ and the set $F(P_w)$ is called the \emph{weak Pareto front}.
\end{enumerate}
\label{def:Pareto_opt}
\end{mydef}
From Definition \ref{def:Pareto_opt}, it immediately follows that $P \subset P_w$. Solving problem \eqref{eq:MOP} means in our setting computing one Pareto optimal point. We do not aim to compute the entire Pareto set. This can be done in a consecutive step using globalization strategies (see \cite{Dellnitz2005}). Since Pareto optimality is defined as a global property, the definition cannot directly be used in practice to check whether a given point is Pareto optimal. Fortunately, there are first order optimality conditions that we can use instead. As in the singleobjective case, they are known as the \emph{Karush-Kuhn-Tucker conditions}.

\begin{mydef}
The set $\Delta^m \coloneqq \left\lbrace \theta \in \R^m \,\,:\,\, \theta \ge 0 \,\,\text{and}\,\, \sum_{i=1}^m \theta_i = 1 \right\rbrace$ is the positive unit simplex. A point $x^* \in \H$ satisfies the \emph{Karush-Kuhn-Tucker conditions} if there exists $\theta \in \Delta^m$ such that 
\begin{align}
\label{eq:KKT_cond}
\sum_{i=1}^m \theta_i \nabla f_i(x^*) = 0.    
\end{align}
If $x^*$ satisfies the Karush-Kuhn-Tucker conditions, we call it \emph{Pareto critical}. The set of all Pareto critical points is the \emph{Pareto critical set,} which we denote by $P_c$.
\end{mydef}
Condition \eqref{eq:KKT_cond} is equivalent to $0 \in \co\left(\left\lbrace \nabla f_i(x^*) \,\,:\,\, i = 1, \dots, m \right\rbrace \right)$ which is also known as Fermat's rule. Analogously to the singleobjective setting, criticality of a point is a necessary condition for optimality. In the convex setting, the KKT conditions are also sufficient conditions for weak Pareto optimality and we have the relation $P\subset P_w = P_c$.

\subsection{Convergence analysis and merit functions}
\label{subsec:conv_analysis_moo}

How do we characterize the convergence of function values for multiobjective optimization problems? For singleobjective optimization problems of the form $\min_{x \in \H} f(x)$ with a convex objective function $f:\H \to \R$, we are interested in the rate of convergence of $f(x(t)) - \inf_{x\in\H} f(x)$. For the problem \eqref{eq:MOP} there is no solution yielding the smallest function value for all objective functions. In the image set $F(\H) = \lbrace F(x) : x \in \H \rbrace \subset \R^m$ there is a set of nondominated points, the so-called Pareto front, which is $F(P) \subset F(\H)$. If we want to characterize the rate of convergence of the function values of a trajectory $x$ for a multiobjective optimization problem, this should relate to the distance of $F(x(t))$ to the Pareto front.

A line of research considers so-called merit functions for multiobjcetive optimization problems to characterize the rate of convergence of function values (see \cite{Tanabe2022, Tanabe2022_2, Tanabe2022_3, Chen2000, Yang2002, Liu2009} and further references in \cite{Tanabe2022_3}). A merit function associated with an optimization problem is a function that returns zero at an optimal solution and which is strictly positive otherwise. In this paper we restrict ourselves to the merit function
\reqnomode
\begin{align}
    u_0(x) \coloneqq \sup_{z \in \H} \min_{i = 1,\dots, m} f_i(x) - f_i(z).
    \label{eq:merit_function}
\end{align}
This function is indeed a merit function for multiobjective optimization problems with respect to weak Pareto optimality, as the following theorem states.
\begin{theorem}
Let $u_0(x)$ be defined by \eqref{eq:merit_function}. For all $x \in \H$ it holds that $u_0(x)\ge 0$. Moreover, $x \in \H$ is weakly Pareto optimal for \eqref{eq:MOP} if and only if $u_0(x) = 0$.
\label{thm:u_0_properties}
\end{theorem}
\begin{proof}
    A proof can be found in \cite[Theorem 3.1]{Tanabe2022_3}.
\end{proof}
Additionally, $u_0$ is lower semicontinuous. Therefore, accumulation points of a smooth curve $t \mapsto x(t)$ with $\lim_{t \to + \infty} u_0(x(t)) = 0$ are weakly Pareto optimal. This motivates the usage of $u_0$ as a measure of convergence speed for multiobjective optimization methods. Even if all objective functions are smooth, the function $u_0$ is in general not smooth. This has to be kept in mind when we define Lyapunov type functions for the system \eqref{eq:MAVD} involving $u_0$.

In the analysis laid out in Section \ref{sec:asym_behaviour_traj} we require the following standing assumption on \eqref{eq:MOP}. This assumption describes a condition on the weak Pareto set.
\begin{myassump}
Let $P_w$ be the set of weakly Pareto optimal points for \eqref{eq:MOP}, and define for $\hat{F} \in \R^m$ the lower level set $\mathcal{L}(\hat{F}) \coloneqq \left\lbrace x \in \H \,\, : \,\, F(x) \le \hat{F}\right\rbrace$. Further, define $P_w(\hat{F}) \coloneqq P_w \cap \mathcal{L}(\hat{F})$. We assume that for all $x_0\in \H$ and all $x \in \mathcal{L}(F(x_0))$ there exists $x^* \in P_w(F(x))$ and
\begin{align*}
    R \coloneqq \sup_{F^* \in F(P_w(F(x_0)))} \inf_{x \in F^{-1}(F^*)} \frac{1}{2}\left\lVert x - x_0 \right\rVert^2 < +\infty.
\end{align*}
\label{assump:1}
\end{myassump}

\begin{myremark}
Assumption \ref{assump:1} is satisfied in the following cases.
\begin{enumerate} [(i)]
    \item For singleobjective optimization problems, Assumption \ref{assump:1} is satisfied if the optimization problem has at least one optimal solution. In this setting, for all $x_0 \in \H$ the weak Pareto set $P_w = P_w(F(x_0))$ coincides with the solution set $\argmin_{x \in \H} f(x)$ and $\inf_{x\in P_w} \frac{1}{2}\lVert x - x_0 \rVert^2 < + \infty$ holds.
    \item Assumption \ref{assump:1} is valid, if the level set $\mathcal{L}(F(x_0))$ is bounded. For example, this is the case when for at least one $i \in \lbrace 1, \dots, m \rbrace$ the set $\left\lbrace x \in \H : f_i(x) \le f_i(x_0) \right\rbrace$ is bounded.
\end{enumerate}
\label{rem:well_behaved_mop}
\end{myremark}
\leqnomode

\begin{myexample}
    We discuss Assumption \ref{assump:1} and Remark \ref{rem:well_behaved_mop} by the means of two examples.
    \begin{enumerate}[(i)]
        \item Assumption \ref{assump:1} cannot hold when the weak Pareto set is empty. Since any scalar optimization problem is also a multiobjective optimization problem, we can consider
        \begin{align}
        \tag{\mbox{MOP$_1$}}
            \min_{x \in \R} F(x) \coloneqq \exp(x).
            \label{eq:P_1}
        \end{align}
        For all $x_0 \in \R$ the set $P_w(F(x_0)) \subset P_w = \argmin_{x\in \R} F(x) = \emptyset$ is empty and hence Assumption \ref{assump:1} does hold. 
        \item In the second example the weak Pareto set is nonempty but the supremum defining $R$ is not bounded. Consider the multiobjective optimization problem
        \begin{align}
        \tag{\mbox{MOP$_2$}}
            \min_{x \in \R^2} F(x) \coloneqq \left[ \begin{array}{c}
                x_1^2  \\
                \exp(x_2)
            \end{array} \right],
            \label{eq:P_2}
        \end{align}
        with two objective functions defined on $\R^2$. For \eqref{eq:P_2} the weak Pareto set is $P_w = \{0\} \times \R$. For all $x_0 \in \R^2$ it holds that $P_w(F(x_0)) = \{ 0 \} \times \left(- \infty, (x_0)_2 \right] \neq \emptyset$, but for this problem $R = + \infty$.
        \end{enumerate}
    For the problems \eqref{eq:P_1} and \eqref{eq:P_2} all objective functions have unbounded level sets. As stated in Remark \ref{rem:well_behaved_mop} (ii), a bounded level set of one objective functions is a sufficient condition for Assumption \ref{assump:1} to hold.
\end{myexample}

\reqnomode

For the convergence analysis in Section \ref{sec:asym_behaviour_traj}, we need an additional lemma on the merit function $u_0$. This lemma describes how to retrieve $u_0$ from $\min_{i=1,\dots,m} f_i(x) - f_i(z)$ without taking the supremum over the whole space $\H$. We need this lemma in particular when we apply the supremum to inequalities in order to bound $u_0$.
\begin{mylemma}
    Let $x_0 \in \H$ and $x \in \mathcal{L}(F(x_0))$, then
    \begin{align*}
        \sup_{F^* \in F(P_w(F(x_0)))} \inf_{z \in F^{-1}(F^*)} \min_{i=1,\dots,m} f_i(x) - f_i(z) = \sup_{z \in \H} \min_{i=1,\dots,m} f_i(x) - f_i(z).
    \end{align*}
    \label{lem:sup_inf_u_0}
\end{mylemma}
\begin{proof}
    A proof of this statement is contained in the proof of Theorem 5.2 in \cite{Tanabe2022}.
\end{proof}

\section{Global existence of solutions to (MAVD) in finite dimensions}
\label{sec:global_existence}
In this section, we show that solutions exist for the Cauchy problem related to
\eqref{eq:MAVD}, i.e.,
\leqnomode
\begin{align}
\label{eq:CP}
\tag{CP}
    \left\lvert 
    \begin{array}{l}
        \frac{\alpha}{t}\dot{x}(t) + \proj_{C(x(t)) + \ddot{x}(t)}{(0)} = 0\, \text{ for } t > t_0,  \\
        \\
        x(t_0) = x_0, \quad \dot{x}(t_0) = v_0,
    \end{array}
    \right.
\end{align}
\reqnomode
with initial data $x_0, v_0 \in \H$ and starting time $t_0 > 0$. The differential equation in \eqref{eq:CP} is implicit. Therefore, we cannot use the Cauchy-Lipschitz or the Peano Theorem to prove existence of solutions. To overcome this problem, we show that solutions exist for \eqref{eq:CP} if there exist solutions to a first order differential inclusion. Then, we use an existence theorem for differential inclusions from \cite{Aubin2012}. Using this approach, we do not expect solutions to be twice continuously differentiable but allow solutions to \eqref{eq:CP} to be less smooth. We specify what we mean by a solution to the Cauchy problem \eqref{eq:CP} in Definition \ref{def:sol_CP}. Since the coefficient $\frac{\alpha}{t}$ has a singularity at $t = 0$, we restrict the analysis in this paper to the case $t_0 > 0$. As our argument only works in finite-dimensional Hilbert spaces, we demand $\dim(\H) < +\infty$ in this section. In our context, the set-valued map
\begin{align}
\label{eq:set_valued_G}
    \hspace{-10mm}G: [t_0, + \infty) \times \H \times \H \rightrightarrows \H\times \H, (t,u,v) \mapsto \lbrace v \rbrace \times \left( - \frac{\alpha}{t} v - \argmin_{g \in C(u)} \langle g, -v \rangle \right)
\end{align}
is of interest. As stated above, $C(u) \coloneqq \co \left(\lbrace \nabla f_i(u) \, : \, i = 1,\dots,m \rbrace \right)$. We can show that \eqref{eq:CP} has a solution if the differential inclusion
\leqnomode
\begin{align}
    \label{eq:DI}
    \tag{DI}
    \left\lvert 
    \begin{array}{l}
        (\dot{u}(t), \dot{v}(t)) \in G(t, u(t), v(t))\, \text{ for } t > t_0, \\
        \\
        (u(t_0), v(t_0)) = (u_0, v_0),
    \end{array}
    \right.
\end{align}
\reqnomode
with appropriate initial data $u_0, v_0 \in \H$ and $t_0 > 0$ has a solution.

\subsection{Existence of solutions to (DI)}
\label{subsec:ex_sol_DI}

To show that there exist solutions to \eqref{eq:DI}, we investigate the set-valued map $(t, u,v) \rightrightarrows G(t, u,v)$ defined in \eqref{eq:set_valued_G}. For a more detailed introduction to the basic definitions regarding set-valued maps, the reader is referred to \cite{Aubin2012}.

\begin{myprop}
    The set-valued map $G$ defined in \eqref{eq:set_valued_G} has the following properties:
    \begin{enumerate}[(i)]
        \item For all $(t,u,v) \in [t_0, +\infty) \times \H \times \H$, the set $G(t,u,v) \subset \H \times \H$ is convex and compact.
        \item $G$ is upper semicontinuous.
        \item For $\dim(\H) < + \infty$, the map 
        \begin{align*}
            \phi: [t_0, +\infty) \times \H \times \H \to \H \times \H, \quad (t,u,v) \mapsto \proj_{G(t,u,v)}(0)  
        \end{align*}
        is locally compact.
        \item Assume the gradients $\nabla f_i$ are globally $L$-Lipschitz continuous. Then, there exists $c > 0$ such that for all $(t, u, v) \in [t_0, +\infty) \times \H \times \H$, it holds that
        \begin{align*}
            \sup_{\xi \in G(t, u, v)} \lVert \xi \rVert_{\H \times \H} \le c(1 + \lVert (u,v) \rVert_{\H \times \H}),
        \end{align*}
        where for $(x,y) \in \H \times \H$ we have $\lVert (x,y) \rVert_{\H \times \H} = \sqrt{\lVert x \rVert^2 + \lVert y \rVert^2}$.
    \end{enumerate}
    \label{prop:set_valued_G}
\end{myprop}
\begin{proof}
    The proof is contained in Appendix \ref{sec:appendix_set_valued_G}.
\end{proof}
The following existence theorem from \cite[p.~98, Theorem 3]{Aubin2012} is applicable in our setting.
\begin{theorem}
\label{thm:existence}
Let $\mathcal{X}$ be a real Hilbert space and let $\Omega \subset \R \times \mathcal{X}$ be an open subset containing $(t_0,x_0)$. Let $G$ be an upper semicontinuous map from $\Omega$ into the nonempty closed convex subsets of $\mathcal{X}$. We assume that $(t,x) \mapsto \proj_{G(t,x)}(0)$ is locally compact. Then, there exists $T > t_0$ and an absolutely continuous function $x$ defined on $[t_0, T]$ which is a solution to the differential inclusion
\begin{align}
    \label{eq:general_diff_incl}
    \dot{x}(t) \in G(t, x(t)), \quad x(t_0) = x_0.
\end{align}
\end{theorem}

In the following remark, we want to give a more precise description of the solutions to a differential inclusion \eqref{eq:general_diff_incl} and give more insight into Theorem \ref{thm:existence}. This is particularly important since the main results of this paper are concerned with the asymptotic behaviour of the solutions of a differential inclusion.

\begin{myremark}
Consider the general differential inclusion \eqref{eq:general_diff_incl}. A solutions $x:[t_0, T] \to \mathcal{X}$ given by Theorem \ref{thm:existence} is not differentiable but merely absolutely continuous. Therefore, the notion $\dot{x}(t) \in G(t, x(t))$ may not hold on the entire domain $[t_0, T]$. An absolutely continuous function $x:[t_0, T] \to \mathcal{X}$ is differentiable almost everywhere in $[t_0, T]$. A solution $x$ to \eqref{eq:general_diff_incl} satisfies the inclusion $\dot{x}(t) \in G(t, x(t))$ in every $t$, where the derivative $\dot{x}(t)$ is defined.
In general $\dot{x}$ will not be continuous. But since $x$ is absolutely continuous with values in a Hilbert space (which satisfies the Radon-Nikodym property) $\dot{x}$ is Bochner integrable and $x(t) = x(t_0) + \int_{t_0}^t \dot{x}(s) ds$ (see \cite{Clarkson1936, Diestel1977}).
\end{myremark}

Within the next theorem, we prove existence of solutions to the differential inclusion \eqref{eq:DI}. 

\begin{theorem}
Assume $\H$ is finite-dimensional and that the gradients of the objective function $\nabla f_i$ are globally Lipschitz continuous. Then, for all $(u_0, v_0) \in \H \times \H$ there exists $T > t_0$ and an absolutely continuous function $(u(\cdot), v(\cdot))$ defined on $[t_0, T]$ which is a solution to the differential inclusion \eqref{eq:DI}.
\label{thm:DI_sol_exist_fin}
\end{theorem}
\begin{proof}
The proof follows immediately from Proposition \ref{prop:set_valued_G} which shows that the set-valued map satisfies all conditions required for Theorem \ref{thm:existence}.  
\end{proof}
In the following, we show that under additional conditions on the objective functions $f_i$, there exist solutions defined on $[t_0, +\infty)$.
\begin{theorem}
Assume $\H$ is finite-dimensional and that the gradients of the objective function $\nabla f_i$ are globally Lipschitz continuous. Then, for all $(u_0, v_0) \in \H \times \H$ there exists a function $(u(\cdot), v(\cdot))$ defined on $[t_0, +\infty)$ which is absolutely continuous on $[t_0, T]$ for all $T > t_0$ and which is a solution to the differential inclusion \eqref{eq:DI}.
\label{thm:DI_sol_exist}
\end{theorem}
\begin{proof}
Theorem \ref{thm:DI_sol_exist_fin} guarantees the existence of solutions defined on $[t_0, T)$ for some $T \ge t_0$. Using the domain of definition, we can define a partial order on the set of solutions of \eqref{eq:DI}. Assuming there is no solution defined on $[t_0, +\infty)$, Zorn's Lemma guarantees the existence of a solution $(u(\cdot), v(\cdot)):[t_0, T) \to \H \times \H$ with $T < + \infty$ which can not be extended. Then we construct an extension using Theorem \ref{thm:DI_sol_exist_fin} to arrive at a contradiction. We give a formal proof built on this idea in the following. Define the set of solutions 
\begin{align*}
    \mathfrak{S} \coloneqq \left\lbrace (u, v):[t_0, T) \to \H \times \H : T \in [t_0, +\infty],\, (u, v) \text{ is a solution to \eqref{eq:DI} on } [t_0, T) \right\rbrace.
\end{align*}
  (Note that $T \in [t_0, +\infty]$ allows for the value $+\infty$ for $T$). Theorem \ref{thm:DI_sol_exist_fin} guarantees $\mathfrak{S} \neq \emptyset$. We define a partial order $\preccurlyeq$ on $\mathfrak{S}$ in the following way. Let $(u_1(\cdot), v_1(\cdot)):[t_0, T_1) \to \H \times \H$ and $ (u_2(\cdot), v_2(\cdot)):[t_0, T_2) \to \H \times \H$ be two solutions to \eqref{eq:DI} with $T_1, T_2 \ge t_0$. Then, we define
\begin{align*}
    &(u_1(\cdot), v_1(\cdot)) \preccurlyeq (u_2(\cdot), v_2(\cdot)) \logeq \\
    T_1 \le T_2 \,\text{ and }\, &(u_1(t), v_1(t)) = (u_2(t), v_2(t)) \,\text{ for all }\, t \in [t_0, T_1).
\end{align*}
This partial order is reflexive, transitive and antisymmetric. We show that any nonempty chain in $\mathfrak{S}$ is bounded by an element of $\mathfrak{S}$. Let $\mathfrak{C} \subseteq \mathfrak{S}$ be an arbitrary nonempty chain, i.e., a totally ordered nonempty subset of $\mathfrak{S}$. From the solutions in $\mathfrak{C}$ we define an upper bound of $\mathfrak{C}$ in $\mathfrak{S}$. Define the supremum
\begin{align*}
    T_{\mathfrak{C}} \coloneqq \sup\left\lbrace T \in [t_0, +\infty] : \text{it exists a solution} (u(\cdot), v(\cdot)) \in \mathfrak{C} \text{ defined on } [t_0, T) \right\rbrace
\end{align*}
Then, we define the solution
\begin{align*}
    &(u_{\mathfrak{C}}(\cdot), v_{\mathfrak{C}}(\cdot)) : [t_0, T_{\mathfrak{C}}) \to \H \times \H, t \mapsto (u(t), v(t)),\\
    \text{ with }&(u(\cdot), v(\cdot)) \in \mathfrak{C}\text{ a solution of \eqref{eq:DI} defined on a strict superset of} \left[t_0, t\right).
\end{align*}
By construction, $(u_{\mathfrak{C}}(\cdot), v_{\mathfrak{C}}(\cdot))$ is well-defined. Further, $(u_{\mathfrak{C}}(\cdot), v_{\mathfrak{C}}(\cdot)) \in \mathfrak{S}$, since at any point $t \in [t_0, T_{\mathfrak{C}})$ the function $(u_{\mathfrak{C}}(\cdot), v_{\mathfrak{C}}(\cdot))$ coincides with a solution from $\mathfrak{C}$ in an open neighborhood of $t$. From the definition of $(u_{\mathfrak{C}}(\cdot), v_{\mathfrak{C}}(\cdot))$ it follows that $(u(\cdot), v(\cdot)) \preccurlyeq (u_{\mathfrak{C}}(\cdot), v_{\mathfrak{C}}(\cdot))$ for all $(u(\cdot), v(\cdot)) \in \mathfrak{C}$, and hence we found an upper bound. We are in the position to use Zorn's Lemma, which guarantees the existence of a maximal element in $\mathfrak{S}$. Let $(u(\cdot), v(\cdot))$ be a maximal element in $\mathfrak{S}$. If $(u(\cdot), v(\cdot))$ is defined on $[t_0, + \infty)$, there exists a global solution and the proof is complete. Therefore, we assume that $(u(\cdot), v(\cdot))$ is defined only on $[t_0, T)$ for some finite real number $T \in [t_0, +\infty)$. We will show that we can extend this solution in order to derive a contradiction to the claimed maximality. To achieve this, we show $(u(t), v(t))$ does not blow up in finite time and use this to construct an extension. Define the function
\begin{align*}
    h(t) \coloneqq \left\lVert (u(t), v(t)) - (u(t_0), v(t_0)) \right\rVert_{\H \times \H}.
\end{align*}
We show that $h(t)$ is bounded by a real-valued function defined on $[t_0, T]$. Using the Cauchy-Schwarz inequality, we get
\reqnomode
\begin{align}
\label{eq:derivative_h_sq}
\begin{split}
    \frac{d}{dt} \frac{1}{2} h^2(t) &= \left\langle (\dot{u}(t), \dot{v}(t)), (u(t), v(t)) - (u(t_0), v(t_0)) \right\rangle_{\H \times \H}\\
    &\le \left\lVert (\dot{u}(t), \dot{v}(t)) \right\rVert_{\H \times \H} h(t).
\end{split}
\end{align}
Proposition \ref{prop:set_valued_G} (iii) guarantees the existence of a constant $c > 0$ with
\begin{align}
\label{eq:bound_h_1}
    \lVert (\dot{u}(t), \dot{v}(t)) \rVert_{\H \times \H} \le c(1 + \lVert (u(t), v(t)) \rVert_{\H \times \H}). 
\end{align}
Define $\tilde{c} \coloneqq c \left( 1 + \lVert (u(t_0), v(t_0)) \rVert_{\H \times \H} \right)$. Then, by applying the triangle inequality to \eqref{eq:bound_h_1}, we have
\begin{align}
\label{eq:bound_h_2}
    \lVert (\dot{u}(t), \dot{v}(t)) \rVert_{\H \times \H} \le \tilde{c}\left(1 + \lVert(u(t), v(t)) - (u(t_0), v(t_0)) \rVert_{\H \times \H}\right).
\end{align}
Combining inequalities \eqref{eq:derivative_h_sq} and \eqref{eq:bound_h_2}, we get for almost all $t \in [t_0, T)$
\begin{align}
\label{eq:bound_derivative}
    \frac{d}{dt}\frac{1}{2}h^2(t) \le \tilde{c} \left(1 + h(t) \right) h(t).
\end{align}
Using a Gronwall-type argument (see Lemma A.4 and Lemma A.5 in \cite{Brezis1973}) as in Theorem 3.5 in \cite{Attouch2015}, we conclude from \eqref{eq:bound_derivative} that for all $\varepsilon > 0$ and for almost all $t \in [t_0, T - \varepsilon]$
\begin{align*}
    h(t) \le \tilde{c}T \exp(\tilde{c}T).
\end{align*}
Since this upper bound is independent of $t$ and $\varepsilon$, the function $h$ is uniformly bounded on $[0, T)$. Then, $u(\cdot)$ and $v(\cdot)$ are uniformly bounded on $[0, T)$ and from \eqref{eq:bound_h_1} we deduce the boundedness of $\dot{u}(\cdot)$ and $\dot{v}(\cdot)$. Together with the fact that $\dot{u}(\cdot)$ and $\dot{v}(\cdot)$ are absolutely continuous, this implies the existence of the integrals
\begin{align*}
    u_T \coloneqq u_0 + \int_{t_0}^T \dot{u}(s) ds \in \H \, \text{ and } \, v_T \coloneqq v_0 + \int_{t_0}^T \dot{v}(s) ds \in \H.
\end{align*}
Then, we define the differential inclusion
\begin{align}
\label{eq:DI_continuation}
    \left\lvert 
    \begin{array}{l}
        (\dot{u}(t), \dot{v}(t)) \in G(t, u(t), v(t))\, \text{ for } t > T, \\
        \\
        (u(T), v(T)) = (u_T, v_T),
    \end{array}
    \right.
\end{align}
and use Theorem \ref{thm:DI_sol_exist_fin} to deduce the existence of $(\hat{u}(\cdot), \hat{v}(\cdot)):[T, T + \delta)\to \H \times \H$ a solution to \eqref{eq:DI_continuation} for some $\delta > 0$. In a next step, we define the function
\begin{align*}
    (u^*(\cdot), v^*(\cdot)):[t_0, T+ \delta) \to \H\times \H,\, t \mapsto \left\lbrace \begin{array}{cl}
        (u(t),v(t)) & \text{ for } t \in [t_0, T) \\
        (\hat{u}(t), \hat{v}(t)) & \text{ for } t \in [T, T+\delta)
    \end{array} \right.
\end{align*}
by connection the solutions. Now, $(u(\cdot), v(\cdot)) \neq (u^*(\cdot), v^*(\cdot)) \in \mathfrak{S}$ and $(u(\cdot), v(\cdot)) \preccurlyeq (u^*(\cdot), v^*(\cdot))$ which forms a contradiction, since $(u(\cdot), v(\cdot))$ is a maximal element of $\mathfrak{S}$.
\end{proof}

\subsection{Existence of solutions to (CP)}
\label{subsec:ex_sol_CP}
Using the findings of the previous subsection, we can proceed with the discussion of the Cauchy problem \eqref{eq:CP}. In this subsection, we show that solutions to the differential inclusion \eqref{eq:DI} immediately give solutions to the Cauchy problem \eqref{eq:CP}. Before we retrieve this result, we have to specify what we understand under a solution to \eqref{eq:CP}. Due to the implicit structure of the equation \eqref{eq:MAVD}, we cannot guarantee the existence of a twice continuously differentiable function $x$ that satisfies equation \eqref{eq:MAVD}. We reduce the requirements in the following sense. 

\begin{mydef}
\label{def:sol_CP}
    We call a function $x:[t_0, +\infty) \to \H$ a solution to the Cauchy problem \eqref{eq:CP} if it satisfies the following conditions.
    \begin{enumerate}[(i)]
        \item $x \in C^1([t_0, +\infty))$, i.e., $x$ is continuously differentiable on $[t_0, +\infty)$.
        \item $\dot{x}$ is absolutely continuous on $[t_0, T]$ for all $T \ge t_0$.
        \item There exists a (Bochner) measurable function $\ddot{x}:[t_0, +\infty) \to \H$ with $\dot{x}(t) = \dot{x}(t_0) + \int_{t_0}^t \ddot{x}(s) ds$ for all $t \ge t_0$.
        \item $\dot{x}$ is differentiable almost everywhere and $\frac{d}{dt} \dot{x} (t) = \ddot{x}(t)$ holds for almost all $t \in [t_0, +\infty)$.
        \item $\frac{\alpha}{t}\dot{x}(t) + \proj_{C(x(t)) + \ddot{x}(t)}(0) = 0$ holds for almost all $t \in [t_0, +\infty)$.
        \item $x(t_0) = x_0$ and $\dot{x}(t_0) = v_0$ hold.
    \end{enumerate}
\end{mydef}

\begin{myremark}
    Conditions (iii) and (iv) are merely consequences of (ii) (see \cite{Clarkson1936, Diestel1977}), since $\dot{x}$ is absolutely continuous on every compact interval $[t_0, T]$ with values in a Hilbert space (which satisfies the Radon-Nikodym property). The (Bochner) measurability of $\ddot{x}$ will be of importance in the analysis of the trajectories.
\end{myremark}

To show that solutions of \eqref{eq:DI} give solutions $x$ which satisfy part $(v)$ of Definition \ref{def:sol_CP}, we need the following auxiliary lemma.

\begin{mylemma}
Let $\H$ be a real Hilbert space, $C \subset \H$ a convex and compact set and $\eta \in \H$ a fixed vector. Then $\xi \in \eta - \argmin_{\mu \in C} \langle \mu, \eta \rangle$ if and only if $\eta = \proj_{C + \xi}(0)$.
\label{lem:proj_1}
\end{mylemma}
\begin{proof}
Let $\xi \in \eta - \argmin_{\mu \in C} \langle \mu, \eta \rangle$. Then, $\eta - \xi \in \argmin_{\mu \in C} \langle \mu, \eta \rangle$ and hence
\begin{align*}
    \langle \eta - \xi , \eta \rangle \le \langle \mu , \eta \rangle \text{ for all } \mu \in C.
\end{align*}
This is equivalent to 
\begin{align*}
    0 \le \langle \mu + \xi - \eta , \eta \rangle \text{ for all } \mu \in C.
\end{align*}
Since $\eta \in C + \xi$ this is equivalent to $\eta = \proj_{C + \xi}(0)$. The other implication follows analogously.
\end{proof}
\begin{theorem}
Let $x_0, v_0 \in \H$ and $t_0 > 0$. Assume $(u(\cdot), v(\cdot)): [t_0, +\infty) \to \H \times \H$ is a solution to \eqref{eq:DI} with $(u(t_0), v(t_0)) = (x_0, v_0)$. Then, it follows that $x(t) \coloneqq u(t)$ satisfies the differential equation
\begin{align*}
    \frac{\alpha}{t} \dot{x}(t) + \proj_{C(x(t)) + \ddot{x}(t)}(0)= 0,
\end{align*}
for almost all $t \in [t_0, +\infty)$ and $x(t_0) = x_0$, $\dot{x}(t_0) = v_0$.
\label{thm:DI_sol_is_CP_sol}
\end{theorem}
\begin{proof}
Since $(u(\cdot), v(\cdot))$ is a solution to \eqref{eq:DI}, the relations
\begin{align}
\label{eq:from_DI_to_CP}
\begin{split} 
    \dot{u}(t) &= v(t) \text{ and}\\
    \dot{v}(t) &\in - \frac{\alpha}{t} v(t) - \argmin_{g \in C(u(t))} \langle g, -v(t)\rangle,
\end{split}
\end{align}
hold for almost all $t \in [t_0, +\infty)$. Using $\frac{\alpha}{t} > 0$, we can write the second line as $\dot{v}(t) \in - \frac{\alpha}{t} v(t) - \argmin_{g \in C(u(t))} \langle g, - \frac{\alpha}{t} v(t) \rangle$. Using Lemma \ref{lem:proj_1} with $\eta = -\frac{\alpha}{t}v(t)$, $C = C(u(t))$ and $\xi = \dot{v}(t)$, the second line in \eqref{eq:from_DI_to_CP} gives
\begin{align*}
    - \frac{\alpha}{t} v(t) = \proj_{C(u(t)) + \dot{v}(t)}(0).
\end{align*}
Rewriting this system using $x(t) = u(t)$, $\dot{x}(t) = \dot{u}(t) = v(t)$ and $\ddot{x}(t) = \dot{v}(t)$ and verifying the initial conditions $x(t_0) = u(t_0) = x_0$ and $\dot{x}(t_0) = v(t_0) = v_0$ yields the desired result. 
\end{proof}
Finally, we can state the full existence theorem for the Cauchy problem \eqref{eq:CP}. 
\begin{theorem}
Assume $\H$ is finite-dimensional and the gradients of the objective function $\nabla f_i$ are globally Lipschitz continuous. Then, for all $x_0, v_0 \in \H$, there exists a function $x$ which is a solution to the Cauchy problem \eqref{eq:CP} in the sense of Definition \ref{def:sol_CP}.
\label{thm:ex_sol_CP}
\end{theorem}
\begin{proof}
The proof follows immediately combining Theorem \ref{thm:DI_sol_exist} and Theorem \ref{thm:DI_sol_is_CP_sol}. 
\end{proof}
\begin{myremark}
In Theorem \ref{thm:ex_sol_CP}, we assume that the gradients $\nabla f_i$ of the objective functions are globally Lipschitz continuous. One can relax this condition and only require the gradients to be Lipschitz continuous on bounded sets if we can guarantee that the solutions remain bounded. This holds for example if one of the objective functions $f_i$ has bounded level sets.
\end{myremark}

\section{Asymptotic behaviour of trajectories}
\label{sec:asym_behaviour_traj}
This section contains the main results regarding the asymptotic properties of the solutions to \eqref{eq:CP}. We show that for $\alpha > 0$ the trajectories $x$ of \eqref{eq:MAVD} minimize the function values. In Theorem \ref{thm:conv_W}, we show that $u_0(x(t)) + \frac{1}{2}\lVert \dot{x}(t) \rVert^2 \to 0$ as $t \to + \infty$ holds in this setting. Hence, every weak accumulation point of $x$ is weakly Pareto optimal. For $\alpha \ge 3$, we prove fast convergence for the function values with rate $u_0(x(t)) = \mathcal{O}(t^{-2})$ as $t \to + \infty$, as we show in Theorem \ref{thm:u_0_x(t)_=_O(t^-2)}. In Theorem \ref{thm:weak_conv}, we prove that for $\alpha > 3$ the trajectories $x$ of \eqref{eq:MAVD} converge weakly to weakly Pareto optimal points using Opial's lemma.

\subsection{Preliminary remarks and estimations}
\label{sec:pre_remarks}
Throughout this subsection, we fix a solution $x:[t_0, +\infty) \to \H$ to \eqref{eq:CP} in the sense of Definition \ref{def:sol_CP} with initial velocity $\dot{x}(t_0) = 0$. Setting the initial velocity to zero has the advantage that the trajectories $x$ remain in the level set $\mathcal{L}(F(x_0))$ as stated in Corollary \ref{cor:descent_wrt_start_value}. Hence, if the level set $\mathcal{L}(F(x_0))$ is bounded, the solution $x$ remains bounded. 
\begin{myprop}
    Let $x:[t_0, +\infty) \to \H$ be a solution to \eqref{eq:CP}. For $i = 1,\dots,m$, define the global energy
    \begin{align}
    \label{eq:def_W_i}
        \W_i:[t_0, +\infty) \to \R, \quad t\mapsto f_i(x(t)) + \frac{1}{2}\lVert \dot{x}(t) \rVert^2.
    \end{align}
    Then, for all $i = 1,\dots, m$ and almost all $t \in [t_0,+\infty)$, it holds that $\frac{d}{dt}\W_i(t) \le -\frac{\alpha}{t}\lVert \dot{x}(t) \rVert^2$. Hence, $\W_i$ is nonincreasing, and $\W_i^{\infty} = \lim_{t \to +\infty} \W_i(t)$ exists in $\R \cup \{-\infty\}$. If $f_i$ is bounded from below, then $\W_i^{\infty} \in \R$.
\label{prop:dissipative_const}
\end{myprop}
\begin{proof}
    The function $\W_i$ is differentiable almost everywhere in $[t_0, +\infty)$ with derivative
    \begin{align}
        \label{eq:W_i_conv_1}
        \frac{d}{dt} \W_i(t) = \frac{d}{dt} \left[f_i(x(t)) +  \frac{1}{2}\lVert \dot{x}(t) \rVert^2 \right] = \left\langle \nabla f_i(x(t)) , \dot{x}(t) \right\rangle + \left\langle \dot{x}(t) , \ddot{x}(t) \right\rangle. 
    \end{align}
    Using the variational representation of $-\frac{\alpha}{t}\dot{x}(t) = \proj_{C(x(t)) + \ddot{x}(t)}(0)$ and the fact that $\nabla f_i(x(t)) \in C(x(t))$, we get for all $i = 1,\dots, m$
    \begin{align*}
        \left\langle \ddot{x}(t) + \frac{\alpha}{t}\dot{x}(t) + \nabla f_i(x(t)), \dot{x}(t) \right\rangle \le 0,
    \end{align*}
    and hence,
    \begin{align}
        \label{eq:W_i_conv_2}
        \left\langle \nabla f_i(x(t)), \dot{x}(t) \right\rangle + \left\langle \ddot{x}(t), \dot{x}(t) \right\rangle \le - \frac{\alpha}{t} \lVert \dot{x}(t) \rVert^2.
    \end{align}
    Combining \eqref{eq:W_i_conv_1} and \eqref{eq:W_i_conv_2} gives the desired results.
\end{proof}

Due to the inertial effects in \eqref{eq:MAVD}, there is in general no monotone descent for the objective values along the trajectories. The following corollary guarantees that the function values along the trajectories are bounded from above by the initial function values given $\dot{x}(t_0) = 0$.
\begin{mycorollary}
     Let $x:[t_0, +\infty) \to \H$ be a solution to \eqref{eq:CP} with $\dot{x}(t_0) = 0$. For all $i = 1,\dots,m$ and all $t \in [t_0, +\infty)$, it holds that
     \begin{align*}
         f_i(x(t)) \le f_i(x_0),
     \end{align*}
     i.e., $x(t) \in \mathcal{L}(F(x_0))$ for all $t \ge t_0$.
     \label{cor:descent_wrt_start_value}
\end{mycorollary}
\begin{proof}
    From Proposition \ref{prop:dissipative_const}, we follow for all $t \in [t_0, +\infty)$
    \begin{align*}
        f_i(x_0) = \W_i(t_0) \ge \W_i(t) = f_i(x(t)) + \frac{1}{2} \lVert \dot{x}(t) \rVert^2 \ge f_i(x(t)).
    \end{align*}
\end{proof}

In the following proofs, we need the weights $\theta(t) \in \Delta^m$ which are implicitly given by
\begin{align}
\label{eq:proj_convex_combination}
    -\frac{\alpha}{t} \dot{x}(t) = \proj_{C(x(t)) + \ddot{x}(t)}(0) = \sum_{i=1}^m \theta_i(t) \nabla f_i(x(t)) + \ddot{x}(t),
\end{align}
    for $t \in [t_0, +\infty)$. In the proofs of the following lemmas, we take the integral over the weights $\theta(t)$. Therefore, we have to guarantee that we can find a measurable selection $t \mapsto \theta(t) \in \Delta^m$ satisfying \eqref{eq:proj_convex_combination}.
    \begin{mylemma}
        Let $x$ be a solution to \eqref{eq:CP}. Then, there exists a measurable function
        \begin{align*}
            \theta: [t_0, +\infty) \to \Delta^m, \quad t \mapsto \theta(t),
        \end{align*}
        with
        \begin{align*}
            \proj_{C(x(t)) + \ddot{x}(t)}(0) = \sum_{i=1}^m \theta_i(t) \nabla f_i(x(t)) + \ddot{x}(t),
        \end{align*}
        for all $t \in [t_0, +\infty)$.
        \label{lem:theta_measurable}
    \end{mylemma}
    \begin{proof}
        Our proof is based on the proof of Proposition 4.6 in \cite{Attouch2015}. Rewrite $\theta(t)$ as a solution to the problem
        \begin{align*}
            \theta(t) \in \argmin_{\theta \in \Delta^m} j(t, \theta), \text{ with } j(t, \theta) \coloneqq \frac{1}{2}\left\lVert \sum_{i = 1}^m \theta_i \nabla f_i(x(t)) + \ddot{x}(t) \right\rVert^2.
        \end{align*}
        We show that $j$ is a Carathéodory integrand. Then, the proof follows from Theorem 14.37 in \cite{Rockafellar2009}, which guarantees the existence of a measurable selection $\theta : [t_0, +\infty) \to \Delta^m, t \mapsto \theta(t) \in \argmin_{\theta \in \Delta^m} j(t, \theta)$. For all $t \in [t_0, +\infty)$, the function $\theta \mapsto j(t, \theta)$ is continuous. By Theorem \ref{thm:ex_sol_CP}, $x$ is a solution to \eqref{eq:CP} in the sense of Definition \ref{def:sol_CP}. This means $\ddot{x}$ is (Bochner) integrable on every interval $[t_0, T]$ and therefore (Bochner) measurable. Then, for all $\theta \in \Delta^m$ the function $t \mapsto j(\theta, t)$ is measurable as a composition of a measurable and a continuous function. This implies that $j$ is indeed a Carathéodory integrand which completes the proof.
    \end{proof}
    In the following, whenever we write $\theta$ we mean the measurable function given by Lemma \ref{lem:theta_measurable}.


\begin{mylemma}
    For $z \in \H$, define $h_z:[t_0, +\infty) \to \R$ by
    \begin{align}
    \label{eq:def_h_z}
        h_z(t) \coloneqq \frac{1}{2}\lVert x(t) - z \rVert^2.
    \end{align}
    For almost all $t \in [t_0, +\infty)$, it holds that
    \begin{align}
        \ddot{h}_z(t) + \frac{\alpha}{t} \dot{h}_z(t) + \sum_{i=1}^m \theta_i(t) (f_i(x(t)) - f_i(z)) \le \lVert \dot{x}(t) \rVert^2.
    \label{eq:ineq_h_z}
    \end{align}
    \label{lem:diff_ineq_h_z}
\end{mylemma}
\begin{proof}
By the chain rule, we have for almost all $t \in [t_0, +\infty)$
\begin{align*}
    \dot{h}_z(t) = \langle x(t) - z, \dot{x}(t) \rangle \quad \text{and} \quad \ddot{h}_z(t) = \langle x(t) - z, \ddot{x}(t) \rangle + \lVert \dot{x}(t) \rVert^2. 
\end{align*}
We combine these expressions with \eqref{eq:proj_convex_combination} to get
\begin{align}
    \begin{split}
        \ddot{h}_z(t) + \frac{\alpha}{t} \dot{h}_z(t) = & \lVert \dot{x}(t) \rVert^2 + \left\langle x(t) - z, - \sum_{i = 1}^m \theta_i(t) \nabla f_i(x(t)) \right\rangle.
    \end{split}
    \label{eq:diff_ineq_h_z_1}
\end{align}
The objective functions $f_i$ are convex and hence $\langle x(t) - z , \nabla f_i(x(t)) \rangle \ge f_i(x(t)) - f_i(z)$. Using this inequality in \eqref{eq:diff_ineq_h_z_1} gives
\begin{align*}
    \ddot{h}_z(t) + \frac{\alpha}{t} \dot{h}_z(t) + \sum_{i=1}^m \theta_i(t) (f_i(x(t)) - f_i(z)) \le \lVert \dot{x}(t) \rVert^2.
\end{align*} 
\end{proof}
Using this lemma, we derive the following relation between $h_z$ and $\W_i$.
\begin{mylemma}
Take $z\in\H$ and let $\W_i$ and $h_z$ be defined by \eqref{eq:def_W_i} and \eqref{eq:def_h_z}, respectively. Then, for all $t \in [t_0, +\infty)$, it holds that
\begin{align*}
    \int_{t_0}^t \frac{1}{s}\sum_{i=1}^m \theta_i(s)\left(\W_i(s) - f_i(z) \right) ds + \frac{3}{2\alpha} \sum_{i=1}^m \theta_i(t)\left(\W_i(t) - f_i(z) \right) \le C_z - \frac{1}{t}\dot{h}_z(t),
\end{align*}
with $C_z \coloneqq (\alpha + 1) \frac{1}{t_0^2}h_z(t_0) + \frac{3}{2\alpha}\max_{i = 1,\dots, m} \left(f_i(x_0) - f_i(z)\right)$.
\label{lem:W_i_h_z}
\end{mylemma}
\begin{proof}
Adding $\frac{1}{2}\lVert \dot{x}(t) \rVert^2$ to inequality \eqref{eq:ineq_h_z} and dividing by $t$, we get for almost all $t \in [t_0, +\infty)$
\begin{align}
    \frac{1}{t}\ddot{h}_z(t) + \frac{\alpha}{t^2} \dot{h}_z(t) + \frac{1}{t}\sum_{i=1}^m \theta_i(t) (\W_i(t) - f_i(z)) \le \frac{3}{2t}\lVert \dot{x}(t) \rVert^2.
    \label{eq:lem:W_i_h_z_proof_1}
\end{align}
We reorder the terms in inequality \eqref{eq:lem:W_i_h_z_proof_1} and integrate from $t_0$ to $t > t_0$, to obtain
\begin{align*}
    & \int_{t_0}^t\frac{1}{s}\sum_{i=1}^m \theta_i(s) (\W_i(s) - f_i(z)) ds \\
    \le & - \int_{t_0}^t \left( \frac{1}{s}\ddot{h}_z(s) + \frac{\alpha}{s^2} \dot{h}_z(s) \right) ds + \int_{t_0}^t \frac{3}{2s}\lVert \dot{x}(s) \rVert^2 ds.
\end{align*}
Integration by parts on the first integral on the right-hand side and using $\dot{h}_z(t_0) = 0$ gives
\begin{align}
\begin{split}
    \le - \frac{1}{t}\dot{h}_z(t) - (\alpha + 1) \int_{t_0}^t \frac{1}{s^2}\dot{h}_z(s) ds  + \int_{t_0}^t \frac{3}{2s}\lVert \dot{x}(s) \rVert^2 ds.
\end{split}
    \label{eq:lem:W_i_h_z_proof_4}
\end{align}
By Proposition \ref{prop:dissipative_const}, we have
\begin{align}
    \int_{t_0}^t \frac{3}{2s}\lVert \dot{x}(s) \rVert^2 ds \le \frac{3}{2\alpha}\sum_{i=1}^m \theta_i(t)\left(\W_i(t_0) - W_i(t)\right).
    \label{eq:lem:W_i_h_z_proof_5}
\end{align}
Applying inequality \eqref{eq:lem:W_i_h_z_proof_5} to \eqref{eq:lem:W_i_h_z_proof_4} and using $\W_i(t_0) = f_i(x_0)$ yields
\begin{align}
\begin{split}
    & \int_{t_0}^t\frac{1}{s}\sum_{i=1}^m \theta_i(s) (\W_i(s) - f_i(z)) ds \\
    \le & - \frac{1}{t}\dot{h}_z(t) - (\alpha + 1) \int_{t_0}^t \frac{1}{s^2}\dot{h}_z(s) ds  + \frac{3}{2\alpha}\sum_{i=1}^m \theta_i(t)\left(f_i(x_0) - \W_i(t)\right).
\end{split}
\label{eq:lem:W_i_h_z_proof_6}
\end{align}
Using integration by parts one more time gives
\begin{align}
    \hspace{-5mm}\int_{t_0}^t \frac{1}{s^2}\dot{h}_z(s) ds = \frac{1}{t^2}h_z(t) - \frac{1}{t_0^2}h_z(t_0) + \int_{t_0}^t \frac{2}{s^3}h_z(s) ds \ge -\frac{1}{t_0^2} h_z(t_0).
\label{eq:lem:W_i_h_z_proof_7}
\end{align}
Combining \eqref{eq:lem:W_i_h_z_proof_6} and \eqref{eq:lem:W_i_h_z_proof_7}, we derive
\begin{align*}
    & \int_{t_0}^t\frac{1}{s}\sum_{i=1}^m \theta_i(s) (\W_i(s) - f_i(z)) ds \\
    \le &  - \frac{1}{t}\dot{h}_z(t) + (\alpha + 1) \frac{1}{t_0^2}h_z(t_0)  + \frac{3}{2\alpha}\sum_{i=1}^m \theta_i(t)\left(f_i(x_0) - \W_i(t)\right) \\
    \le &  - \frac{1}{t}\dot{h}_z(t) + (\alpha + 1) \frac{1}{t_0^2}h_z(t_0)  + \frac{3}{2\alpha}\sum_{i=1}^m \theta_i(t)\left(f_i(x_0) - f_i(z)\right) \\
    + & \frac{3}{2\alpha}\sum_{i=1}^m \theta_i(t)\left(f_i(z) - \W_i(t)\right) \\
    \le & C_z - \frac{1}{t}\dot{h}_z(t) - \frac{3}{2\alpha} \sum_{i=1}^m \theta_i(t) \left(\W_i(t) - f_i(z)\right),
\end{align*}
with
\begin{align}
\label{eq:def_C_z}
    C_z \coloneqq (\alpha + 1) \frac{1}{t_0^2}h_z(t_0) + \frac{3}{2\alpha}\max_{i = 1,\dots, m} \left(f_i(x_0) - f_i(z)\right),
\end{align}
which completes the proof. 
\end{proof}


\begin{mylemma}
Take $z\in\H$ and let $\W_i, h_z$ and $C_z$ be defined by \eqref{eq:def_W_i}, \eqref{eq:def_h_z} and \eqref{eq:def_C_z}, respectively. Then, for all $\tau > t_0$ it holds that
\begin{align*}
    \min_{i = 1,\dots,m} \left(\W_i(\tau) - f_i(z) \right) \left[ \tau \ln \tau + A \tau + B \right] \le C_z( \tau - t_0) + \frac{h_z(t_0)}{t_0},
\end{align*}
with constants $A, B \in \R$ which are independent of $z$.
\label{lem:W_i_h_z_2}
\end{mylemma}
\begin{proof}
        Let $z \in \H$ and $\tau \ge t > t_0$. Proposition \ref{prop:dissipative_const} states that the functions $\W_i$ are nonincreasing. Therefore, we have for all $s \in [t_0, t]$, that $\W_i(\tau) - f_i(z) \le \W_i(s) - f_i(z)$ and hence
    \begin{align}
        \begin{split}
        & \min_{i = 1,\dots,m} \left(\W_i(\tau) - f_i(z) \right) \int_{t_0}^t \frac{1}{s} ds + \frac{3}{2\alpha}\min_{i = 1,\dots,m} \left(\W_i(\tau) - f_i(z) \right) \\
        \le & \int_{t_0}^t \frac{1}{s} \min_{i = 1,\dots,m} \left(\W_i(s) - f_i(z) \right) ds + \frac{3}{2\alpha}\min_{i = 1,\dots,m} \left(\W_i(t) - f_i(z) \right).
        \end{split}
        \label{eq:minimizing_prop_1}
    \end{align}
    Using Lemma \ref{lem:W_i_h_z}, we get
    \begin{align}
        \begin{split}
        & \int_{t_0}^t \frac{1}{s} \min_{i = 1,\dots,m} \left(\W_i(s) - f_i(z) \right) ds + \frac{3}{2\alpha}\min_{i = 1,\dots,m} \left(\W_i(t) - f_i(z) \right)\\
        \le & \int_{t_0}^t \frac{1}{s} \sum_{i=1}^m \theta_i(s) \left(\W_i(s) - f_i(z) \right) ds + \frac{3}{2\alpha}\sum_{i=1}^m \theta_i(t)\left(\W_i(t) - f_i(z) \right)\\
        \le & C_z - \frac{1}{t}\dot{h}_z(t).
        \end{split}
        \label{eq:minimizing_prop_2}
    \end{align}
    Inequalities \eqref{eq:minimizing_prop_1} and \eqref{eq:minimizing_prop_2} together give
    \begin{align}
        \min_{i = 1,\dots,m} \left(\W_i(\tau) - f_i(z) \right) \left[ \ln t - \ln t_0 + \frac{3}{2\alpha}\right] \le C_z - \frac{1}{t}\dot{h}_z(t).
        \label{eq:minimizing_prop_3}
    \end{align}
    Integrating inequality \eqref{eq:minimizing_prop_3} from $t = t_0$ to $t = \tau$, we have
    \begin{align}
    \begin{split}
        &\hspace{-15mm} \min_{i = 1,\dots,m} \left(\W_i(\tau) - f_i(z) \right) \left[ \tau \ln \tau - \tau - t_0 \ln t_0 + t_0 + \left(\frac{3}{2\alpha} - \ln t_0 \right)(\tau - t_0) \right] \\
        &\hspace{-20mm}\le C_z( \tau - t_0) - \int_{t_0}^{\tau} \frac{1}{t} \dot{h}_z(t) dt.
        \label{eq:minimizing_prop_4}
    \end{split}
    \end{align}
    Integration by parts yields
    \begin{align}
        \int_{t_0}^{\tau} \frac{1}{t} \dot{h}_z(t) dt = \frac{h_z(\tau)}{\tau} - \frac{h_z(t_0)}{t_0} + \int_{t_0}^{\tau} \frac{h_z(t)}{t^2} dt \ge - \frac{h_z(t_0)}{t_0}.
        \label{eq:minimizing_prop_5}
    \end{align}
    Using inequality \eqref{eq:minimizing_prop_5} in \eqref{eq:minimizing_prop_4}, we write
    \begin{align*}
        & \min_{i = 1,\dots,m} \left(\W_i(\tau) - f_i(z) \right) \left[ \tau \ln \tau - \tau - t_0 \ln t_0 + t_0 + \left(\frac{3}{2\alpha} - \ln t_0 \right)(\tau - t_0) \right]\\
        &\hspace{-5mm}\le C_z( \tau - t_0) + \frac{h_z(t_0)}{t_0}.
    \end{align*}
    Introducing suitable constants $A, B\in \R$, this gives the desired result.
\end{proof}


The next theorem is the main result of this subsection. Theorem \ref{thm:conv_W} states that the function values of the trajectories $F(x(t)) \in \R^m$ converge to elements of the Pareto front. In addition, Theorem \ref{thm:conv_W} states that every weak limit point of the trajectory $x$ is weakly Pareto optimal. This is important for proving the weak convergence of the trajectories to weakly Pareto optimal points in Subsection \ref{subsec:weak_conv_traj}. 
\begin{theorem}
    Let $\alpha > 0$ and suppose $x:[t_0, +\infty) \to \H$ is a solution of \eqref{eq:CP}. Define the global energy
    \begin{align*}
        \W(t) = u_0(x(t)) + \frac{1}{2}\lVert \dot{x}(t) \rVert^2 = \sup_{z \in \H} \min_{i=1,\dots,m} \left(f_i(x(t)) - f_i(z) \right) + \frac{1}{2}\lVert \dot{x}(t) \rVert^2.
    \end{align*}
    Assume the functions $f_i$ are bounded from below and Assumption \ref{assump:1} holds. Then, $\lim_{t \to + \infty} \W(t) = 0$. Hence, $\lim_{t \to + \infty} u_0(x(t)) = 0$ and by Theorem \ref{thm:u_0_properties} every weak limit point of $x$ is Pareto critical.
    \label{thm:conv_W}
\end{theorem}

\begin{proof}
    Lemma \ref{lem:W_i_h_z_2} states
    \begin{align}
        \min_{i = 1,\dots,m} \left(\W_i(\tau) - f_i(z) \right) \left[ \tau \ln \tau + A \tau + B \right] \le C_z( \tau - t_0) + \frac{h_z(t_0)}{t_0},
        \label{eq:proof_conv_W_1}
    \end{align}
    for all $\tau > t_0$. We cannot directly take the supremum on both sides since $C_z$ might be unbounded w.r.t. $z \in \H$.
    For $z \in \mathcal{L}(F(x_0))$ we have $\max_{i=1, \dots, m}\left( f_i(x_0) - f_i(z) \right) \le \max_{i=1, \dots, m}\left(f_i(x_0) - \inf_{z \in \H} f_i(z) \right) \eqqcolon M$. Since all $f_i$ are bounded from below by assumption, we have $M < +\infty$.
    Fix $F^* \in F(P_w(F(x_0)))$. Using the definition of $C_z$ given in \eqref{eq:def_C_z}, we get for all $z \in F^{-1}(F^*)$ 
    \begin{align}
    \begin{split}
        C_z(\tau - t_0) + \frac{h_z(t_0)}{t_0} \le & \left((\alpha + 1)\frac{1}{t_0^2}h_z(t_0) + \frac{3}{2 \alpha}M\right)(\tau - t_0) + \frac{h_z(t_0)}{t_0} \\
        = & \left( \frac{\alpha + 1}{t_0^2}(\tau - t_0) + \frac{1}{t_0}\right) h_z(t_0) + \frac{3}{2 \alpha} M(\tau - t_0).
    \end{split}
    \label{eq:proof_conv_W_2}
    \end{align}
    By Assumption \ref{assump:1} $\sup_{F^* \in F(P_w(F(x_0)))} \inf_{z \in F^{-1}(F^*)} h_z(t_0) = R < + \infty$. Applying this infimum and supremum to \eqref{eq:proof_conv_W_2} we have
    \begin{align}
    \begin{split}
        & \sup_{F^* \in F(P_w(F(x_0)))} \inf_{z \in F^{-1}(F^*)} \left( C_z(\tau - t_0) + \frac{h_z(t_0)}{t_0} \right)\\
        \le & \left(\frac{(\alpha + 1)R}{t_0^2} + \frac{3}{2 \alpha}M \right)(\tau - t_0) + \frac{R}{t_0}. 
        \label{eq:proof_conv_W_3}
    \end{split}
    \end{align}
    Combining Lemma \ref{lem:sup_inf_u_0} with \eqref{eq:proof_conv_W_1} and \eqref{eq:proof_conv_W_3}, we get for all $\tau > t_0$
    \begin{align}
    \label{eq:log_ineq_W_tau}
        \W(\tau) \left[ \tau \ln \tau + A \tau + B \right] \le C\tau + D,
    \end{align}
    with $A, B, D \in \R$ and $C > 0$. Since $\W(t)$ is nonnegative $\lim_{t \to + \infty} \W(t) = 0$ holds.
\end{proof}

\begin{myremark}
    From the proof of Theorem \ref{thm:conv_W}, we can deduce a slightly stronger result. There is not only convergence $\lim_{t \to +\infty} \W(t) = 0$, but from inequality \eqref{eq:log_ineq_W_tau} we get convergence of order
    \begin{align*}
        \W(t) = \mathcal{O}\left(\frac{1}{\ln t}\right), \text{ for } t \to + \infty.
    \end{align*}
    Since this is a rather slow rate of convergence, which is not used in the following proofs, we do not point it out in Theorem \ref{thm:conv_W}.
\end{myremark}

We can derive some additional facts on the function values $f_i(x(t))$ along the trajectories from Theorem \ref{thm:conv_W}.

\begin{theorem}
    Let $x:[t_0, +\infty) \to \H$ be a solution to \eqref{eq:CP} and assume all assumptions of Theorem \ref{thm:conv_W} hold. Then, for all $i = 1,\dots, m$
    \begin{align*}
        \lim_{t \to \infty} f_i(x(t)) = f_i^{\infty} \in \R
    \end{align*}
    exists, and further $f_i^{\infty} = \W_i^{\infty}$.
    \label{thm:conv_f_i}
\end{theorem}
\begin{proof}
    Theorem \ref{thm:conv_W} states $\lim_{t \to \infty} \W(t) = 0$. By definition, we have $\W(t) = u_0(x(t)) + \frac{1}{2}\lVert \dot{x}(t) \rVert^2$. Theorem \ref{thm:u_0_properties} guarantees $u_0(x(t)) \ge 0$ for all $t \ge t_0$ and obviously $\frac{1}{2}\lVert \dot{x}(t) \rVert^2 \ge 0$. Then, from $\lim_{t \to + \infty} \W(t) = 0$ it follows that $\lim_{t \to +\infty} \frac{1}{2}\lVert \dot{x}(t) \rVert^2 = 0$. Since $\lim_{t \to +\infty} \W_i(t) = \lim_{t \to + \infty} f_i(x(t)) + \frac{1}{2}\lVert \dot{x}(t) \rVert^2$ exists, it follows that the limit $\lim_{t \to + \infty} f_i(x(t))$ exists, which completes the proof. 
\end{proof}

\subsection{Fast convergence of function values}
\label{subsec:fast_conv_func_val}
In this subsection, we show that solutions of \eqref{eq:CP} have good properties with respect to multiobjective optimization. Along the trajectories of \eqref{eq:CP} the function values converge with order $\mathcal{O}(t^{-2})$ to an optimal value, given $\alpha \ge 3$. This convergence has to be understood in terms of the merit function $u_0$. We prove this result using Lyapunov type energy functions similar to the analysis for the singleobjective case laid out in \cite{Su2014} and \cite{Attouch2018}. To this end, we introduce two important auxiliary functions in Definition \ref{def:energy_functions} and discuss their basic properties in the following lemmas. The main result of this subsection on the convergence of the function values is stated in Theorem \ref{thm:u_0_x(t)_=_O(t^-2)}.
\begin{mydef}
    \label{def:energy_functions}
    Let $\lambda \ge 0, \xi \ge 0$, $z \in \H$ and $x:[t_0, +\infty) \to \H$ be a solution to \eqref{eq:CP}. For $t \ge t_0$  and $i = 1,\dots,m$ define
    \begin{align}
        \hspace{-10mm}\E_{i, \lambda, \xi}(t) = t^2(f_i(x(t)) - f_i(z)) + \frac{1}{2}\lVert \lambda(x(t) - z) + t\dot{x}(t) \rVert^2 + \frac{\xi}{2}\lVert x(t) - z \rVert^2.
        \label{eq:energy_E_i}
    \end{align}
    Using the functions $\E_{i, \lambda, \xi}$, define for $t \ge t_0$
    \begin{align}
        \begin{split}
            &\hspace{-11mm}\E_{\lambda, \xi}(t) = \min_{i=1,\dots,m} \E_{i, \lambda, \xi}(t) \\
            &\hspace{-16mm} = t^2 \min_{i=1,\dots, m} \left( f_i(x(t)) - f_i(z) \right) + \frac{1}{2}\lVert \lambda (x(t) - z ) + t\dot{x}(t) \rVert^2 + \frac{\xi}{2}\lVert x(t) - z \rVert^2.
        \end{split}
        \label{eq:energy_min_E_i}
    \end{align}
\end{mydef}

\begin{mylemma}
    Let $\xi^* = \lambda(\alpha - 1 - \lambda)$ and $\alpha \ge \lambda +1$. Then, for all $i = 1,\dots,m$ and almost all $t \in [t_0, +\infty)$
    \begin{align*}
        \frac{d}{dt} \E_{i, \lambda, \xi^*}(t) \le 2t(f_i(x(t)) - f_i(z)) - t \lambda \min_{i=1,\dots,m}\left( f_i(x(t)) - f_i(z) \right)\\
         + t (\lambda + 1 - \alpha)\lVert \dot{x}(t) \rVert^2.
    \end{align*}
    \label{lem:derivative_E_i_lam_xi}
\end{mylemma}
\begin{proof}
The function $\E_{i, \lambda, \xi^*}(t)$ is differentiable almost everywhere since $f_i$ and $x$ are differentiable and $\dot{x}$ is absolutely continuous. We compute $\frac{d}{dt}\E_{i, \lambda, \xi^*}(t)$ using the chain rule on \eqref{eq:energy_E_i}
\begin{align}
    \begin{split}
    & \frac{d}{dt} \E_{i, \lambda, \xi^*}(t) = 2t(f_i(x(t)) - f_i(z)) + t^2 \left\langle \dot{x}(t) , \nabla f_i(x(t)) + \ddot{x}(t) \right\rangle\\
    & + t \left\langle x(t) - z , \frac{\lambda(\lambda + 1) + \xi^*}{t}\dot{x}(t) +  \lambda \ddot{x}(t) \right\rangle + t(\lambda + 1) \lVert \dot{x}(t) \rVert^2.
    \end{split}
    \label{eq:derivative_E_i_lam_xi_1}
\end{align}
Using Proposition \ref{prop:dissipative_const} on the second summand in \eqref{eq:derivative_E_i_lam_xi_1}, we bound this by
\begin{align}
\label{eq:derivative_E_i_lam_xi_1_0}
\begin{split}
     \le & 2t(f_i(x(t)) - f_i(z)) + t \left\langle x(t) - z , \frac{\lambda(\lambda + 1) + \xi^*}{t}\dot{x}(t)+ \lambda \ddot{x}(t) \right\rangle\\
     & + t(\lambda + 1 - \alpha) \lVert \dot{x}(t) \rVert^2.
\end{split}
\end{align}
We rewrite \eqref{eq:derivative_E_i_lam_xi_1_0} into
\begin{align*}      
     = 2t(f_i(x(t)) - f_i(z)) + t \lambda \left\langle x(t) - z , \frac{\alpha}{t}\dot{x}(t) + \ddot{x}(t) \right\rangle +  t(\lambda + 1 - \alpha) \lVert \dot{x}(t) \rVert^2,
\end{align*}
using $\lambda(\lambda + 1) + \xi^* = \lambda \alpha$. The definition of \eqref{eq:MAVD} together with Lemma \ref{lem:theta_measurable} implies
\begin{align*}
     = 2t(f_i(x(t)) - f_i(z)) - t \lambda \sum_{i=1}^m \theta_i(t) \left\langle x(t) - z , \nabla f_i(x(t)) \right\rangle + t(\lambda + 1 - \alpha) \lVert \dot{x}(t) \rVert^2.
\end{align*}
The objective functions $f_i$ are convex and hence $f_i(z) - f_i(x(t)) \ge \langle \nabla f_i(x(t)), z - x(t) \rangle$ and therefore,
\begin{align*}
     \le 2t(f_i(x(t)) - f_i(z)) - t \lambda \sum_{i=1}^m \theta_i(t)\left( f_i(x(t)) - f_i(z) \right) + t(\lambda + 1 - \alpha) \lVert \dot{x}(t) \rVert^2.
\end{align*}
We bound the convex combination using the minimum to get
\begin{align*}
     \le 2t(f_i(x(t)) - f_i(z)) - t \lambda \min_{i=1,\dots,m}\left( f_i(x(t)) - f_i(z) \right) + t(\lambda + 1 - \alpha) \lVert \dot{x}(t) \rVert^2.
\end{align*}  
\end{proof}

To retrieve a result similar to Lemma \ref{lem:derivative_E_i_lam_xi} for the function $\E_{\lambda, \xi}$ defined in \eqref{eq:energy_min_E_i}, we need an auxiliary lemma which helps us to treat the derivative of $\E_{\lambda, \xi}$.

\begin{mylemma}
    Let $(h_i)_{i=1,\dots,m}$ be a family of continuously differentiable functions $h_i:[t_0, +\infty) \to \R$. Define $h:[t_0, +\infty) \to \R,\quad t \mapsto h(t) \coloneqq \min_{i=1,\dots, m}h_i(t)$. Then, it holds that
    \begin{enumerate}[(i)]
        \item For almost all $t \in [t_0, +\infty)$, the function $h$ is differentiable in $t$. 
        \item For almost all $t \in [t_0, +\infty)$, there exists $i \in \{1,\dots, m\}$ with $h(t) = h_i(t)$ and 
        \begin{align*}
            \frac{d}{dt}h(t) = \frac{d}{dt}h_i(t).    
        \end{align*}
    \end{enumerate}
    \label{lem:diff_a_e}
\end{mylemma}
\begin{proof}

    (i) The functions $h_i$ are continuously differentiable. Therefore, $h$ is locally Lipschitz continuous. Then, by Rademacher's Theorem $h$ is differentiable almost everywhere.
    
    (ii) This statement follows from (i) and the definition of the derivative using the difference quotient.
\end{proof}

\begin{mylemma}
    Let $x:[t_0, +\infty) \to \H$ be a solution to \eqref{eq:CP} and let $z \in \H$. The energy function $\E_{\lambda, \xi^*}$ satisfies the following conditions.
    \begin{enumerate}[(i)]
    \item The function $\E_{\lambda, \xi^*}$ is differentiable in almost all $t \in [t_0, +\infty)$.
    \item For almost all $t \in [t_0, +\infty)$, it holds that
    \begin{align}
        \hspace{-17mm}\frac{d}{dt}\E_{\lambda, \xi^*}(t) \le (2-\lambda)t \min_{i=1,\dots, m} \left( f_i(x(t)) - f_i(z) \right) - (\alpha - \lambda - 1) t \lVert \dot{x}(t) \rVert^2.
        \label{eq:energy_estimation_-1}
    \end{align}
    \item For all $t \in [t_0, +\infty)$, it holds that
    \begin{align*}
        \E_{\lambda, \xi^*}(t) - \E_{\lambda, \xi^*}(t_0) \le (2 - \lambda) & \int_{t_0}^t t \min_{i=1,\dots,m}\left( f_i(x(t)) - f_i(z) \right) dt\\
        + & \int_{t_0}^t t(\lambda + 1 - \alpha) \lVert \dot{x}(t) \rVert^2 dt.
    \end{align*}
    \end{enumerate}
    \label{lem:energy_estimation}
\end{mylemma}
\begin{proof}~\\

    (i)    The functions $t \mapsto f_i(x(t))$ are continuously differentiable for all $i = 1,\dots,m$. Then, by Lemma \ref{lem:diff_a_e} the function $t \mapsto \min_{i=1,\dots,m} \left(f_i(x(t)) - f_i(z) \right)$ is differentiable in $t$ for almost all $t \in [t_0, +\infty)$. Since $x$ is a solution to \eqref{eq:CP} in the sense of Definition \ref{def:sol_CP}, we know that $\lVert \lambda (x(t) - z) + t \dot{x}(t) \rVert^2$ and $\frac{\xi}{2}\lVert x(t) - z \rVert^2$ are differentiable in $t$ for almost all $t \in [t_0, +\infty)$. In total we get that $\E_{\lambda, \xi^*}(t)$ is differentiable in $t$ for almost all $t \in [t_0, +\infty)$.
    
    (ii)    In order to compute the derivative of $\E_{\lambda, \xi^*}(t)$, we need the derivative of $\min_{i=1,\dots,m}\left(f_i(x(t))-f_i(z)\right)$. By Lemma \ref{lem:diff_a_e} for almost all $t \in [t_0, +\infty)$ there exists $j \in \{1, \dots, m\}$ with
    \begin{align}
    \begin{split}
        \frac{d}{dt}\min_{i=1,\dots, m} \left( f_i(x(t)) - f_i(z) \right) &= \frac{d}{dt}\left( f_j(x(t)) - f_j(z) \right), \text{ and } \\
         \min_{i=1,\dots,m} \left(f_i(x(t)) - f_i(z) \right) &= f_j(x(t)) - f_j(z).
    \end{split}
        \label{eq:energy_estimation_1}
    \end{align}
    For the remainder of the proof fix $t \in [t_0, +\infty)$ and $j$ satisfying equation \eqref{eq:energy_estimation_1}. From the first part of \eqref{eq:energy_estimation_1}, we immediately get
    \begin{align}
        \frac{d}{dt}\E_{\lambda, \xi^*}(t) = \frac{d}{dt}\E_{j, \lambda, \xi^*}(t).
        \label{eq:energy_estimation_2}
    \end{align}
    Applying Lemma \ref{lem:derivative_E_i_lam_xi}, we bound \eqref{eq:energy_estimation_2} by 
    \begin{align*}
        \le 2t(f_j(x(t)) - f_j(z)) - t \lambda \min_{i=1,\dots,m}\left( f_i(x(t)) - f_i(z) \right) + t(\lambda + 1 - \alpha) \lVert \dot{x}(t) \rVert^2.
    \end{align*}
    Then, the second equation in \eqref{eq:energy_estimation_1} gives
    \begin{align*}
        = & (2 - \lambda) t \min_{i=1,\dots,m}\left( f_i(x(t)) - f_i(z) \right) + t(\lambda + 1 - \alpha) \lVert \dot{x}(t) \rVert^2.
    \end{align*}
    Statement (iii) follows immediately from (ii) by integrating inequality \eqref{eq:energy_estimation_-1} from $t_0$ to $t$.
\end{proof}
The term $\min_{i=1,\dots, m} \left( f_i(x(t)) - f_i(z) \right)$ will not remain nonnegative in general. Hence, we cannot guarantee that $\E_{\lambda, \xi^*}(t)$ is nonnegative. Therefore, the function $\E_{\lambda, \xi^*}(t)$ is not suitable for convergence analysis and we cannot directly retrieve results on the convergence rates. We are still able to get convergence results using Lemma \ref{lem:sup_inf_u_0}.
\begin{theorem}
    Let $\alpha \ge 3$ and $x:[t_0, + \infty) \to \H$ be a solution to \eqref{eq:CP}. Then
    \begin{align*}
         t^2 u_0(x(t)) \le  t_0^2 u_0(x_0) + 2(\alpha - 1)R + (3 - \alpha) \int_{t_0}^t s \lVert \dot{x}(s) \rVert^2 ds,
    \end{align*}
    and hence $u_0(x(t)) \le \frac{t_0^2 u_0(x_0) + 2(\alpha - 1)R}{t^2}$ for all $t \in [t_0, +\infty)$.
    \label{thm:u_0_x(t)_=_O(t^-2)}
\end{theorem}

\begin{proof}
    We consider the energy function $\E_{\lambda, \xi^*}(t)$ with parameter $\lambda = 2$. From the definition of $\E_{2, \xi^*}(t)$ and part (iii) of Lemma \ref{lem:energy_estimation}, we deduce
    \begin{align*}
         t^2\min_{i=1,\dots,m} \left(f_i(x(t)) - f_i(z) \right) \le \E_{2, \xi^*}(t_0) + (3-\alpha) \int_{t_0}^t s \lVert \dot{x}(s) \rVert^2 ds.
    \end{align*}
    Writing out the definition of $\E_{2, \xi^*}(t_0)$ and using $\lambda = 2$ and $\xi^* = \lambda(\alpha - 1 - \lambda) = 2(\alpha - 3)$, we have 
    \begin{align}
    \begin{split} 
        t^2\min_{i=1,\dots,m} \left(f_i(x(t)) - f_i(z) \right) \le & t_0^2\min_{i=1,\dots,m} \left(f_i(x_0) - f_i(z) \right) + (\alpha - 1)\lVert x_0 - z \rVert^2 \\
        & + (3-\alpha) \int_{t_0}^t s \lVert \dot{x}(s) \rVert^2 ds.
    \end{split}
    \label{eq:u_0_x(t)_=_O(t^-2)_3}
    \end{align}
    We want to apply the supremum and infimum in accordance with Lemma \ref{lem:sup_inf_u_0}. Let $F^* = (f_1^*, \dots, f_m^*) \in F(P_w(F(x_0)))$, then
    \begin{align}
    \begin{split}
        & \inf_{z \in F^{-1}(F^*)}\left[ t_0^2 \min_{i=1,\dots, m} \left(f_i(x_0) - f_i(z)\right) + (\alpha - 1) \lVert x_0 - z \rVert^2 \right] \\
        = & t_0^2 \min_{i=1,\dots, m} \left(f_i(x_0) - f^*_i\right) + (\alpha - 1) \inf_{z \in F^{-1}(F^*)} \lVert x_0 - z \rVert^2.
    \end{split}
    \label{eq:u_0_x(t)_=_O(t^-2)_4}
    \end{align}
    Now, we can apply the supremum to inequality \eqref{eq:u_0_x(t)_=_O(t^-2)_4} and get
    \begin{align}
    \begin{split}
        & \sup_{F^* \in P_w(F(x_0))} \inf_{z \in F^{-1}(F^*)}\left[ t_0^2 \min_{i=1,\dots, m} \left(f_i(x_0) - f_i(z)\right) + (\alpha - 1) \lVert x_0 - z \rVert^2 \right] \\
        \le & t_0^2 \sup_{F^* \in P_w(F(x_0))} \inf_{z \in F^{-1}(F^*)} \min_{i=1,\dots, m} \left(f_i(x_0) - f^*_i\right)\\
        & + (\alpha - 1)  \sup_{F^* \in P_w(F(x_0))} \inf_{z \in F^{-1}(F^*)} \inf_{z \in F^{-1}(F^*)} \lVert x_0 - z \rVert^2.
    \end{split}
    \label{eq:u_0_x(t)_=_O(t^-2)_5}
    \end{align}
    By Assumption \ref{assump:1} and the definition of $u_0(x)$, this is equal to 
    \begin{align}    
        = t_0^2 u_0(x_0) + 2(\alpha - 1)R.
        \label{eq:u_0_x(t)_=_O(t^-2)_6}
    \end{align}
    Now, by applying $\sup_{F^* \in P_w(F(x_0))} \inf_{z \in F^{-1}(F^*)}$ to $t^2 \min_{i = 1,\dots,m}\left(f_i(x(t)) - f_i(z)\right)$ and using \eqref{eq:u_0_x(t)_=_O(t^-2)_3} - \eqref{eq:u_0_x(t)_=_O(t^-2)_6}, we get
    \begin{align*}
        t^2 u_0(x(t)) \le  t_0^2 u_0(x_0) + 2(\alpha - 1)R + (3 - \alpha) \int_{t_0}^t s \lVert \dot{x}(s) \rVert^2 ds.
    \end{align*}
\end{proof}

\begin{mycorollary}
    Let $\alpha > 3$ and $x:[t_0, + \infty) \to \H$ be a solution to \eqref{eq:MAVD}. Then 
    \begin{align*}
        \int_{t_0}^{+\infty} s \lVert \dot{x}(s) \rVert^2 ds < + \infty,
    \end{align*}
    i.e., $(t \mapsto t \lVert \dot{x}(t) \rVert^2) \in L^1([t_0, + \infty))$.
\end{mycorollary}

\subsection{Weak convergence of trajectories}
\label{subsec:weak_conv_traj}
In this subsection, we show that the bounded trajectiories of \eqref{eq:MAVD} converge weakly to weakly Pareto optimal solutions of \eqref{eq:MOP}, given $\alpha > 3$. We prove this in Theorem \ref{thm:weak_conv} using Opial's lemma. Since we need to apply Theorem \ref{thm:conv_W} and Theorem \ref{thm:conv_f_i}, we assume in this subsection that the functions $f_i$ are bounded from below and that Assumption \ref{assump:1} holds. We start this subsection by reciting Opial's lemma (see \cite{Opial1967}) that we need in order to prove weak convergence.
\begin{mylemma}
Let $S \subset \H$ be a nonempty subset of $\H$ and $x:[t_0, +\infty) \to \H$. Assume that $x$ satisfies the following conditions.
\begin{enumerate}[(i)]
    \item Every weak sequential cluster point of $x$ belongs to $S$.
    \item For every $z\in S$, $\lim_{t\to + \infty} \lVert x(t) - z \rVert$ exists.
\end{enumerate}
Then, $x(t)$ converges weakly to an element $x^{\infty} \in S$, as $t \to +\infty$.
\label{lem:opial}
\end{mylemma}
We need the following additional lemma in order to utilize Opial's lemma.
\begin{mylemma}
    Let $t_0 > 0$ and let $h:[t_0, +\infty) \to \R$ be a continuously differentiable function which is bounded from below. Assume
    \begin{align*}
        t \ddot{h}(t) + \alpha \dot{h}(t) \le g(t),
    \end{align*}
    for some $\alpha > 1$ and almost all $t \in [t_0, +\infty)$, where $g \in L^1([t_0, +\infty))$ is a nonnegative function. Then, $\lim_{t \to + \infty} h(t)$ exists.
    \label{lem:differential_inequality}
\end{mylemma}
\begin{proof}
    A proof can be found in \cite[Lemma A.6.]{Attouch2019_2}.
\end{proof}

\begin{theorem}
    Let $\alpha > 3$ and let $x:[t_0, +\infty) \to \H$ be a bounded solution to \eqref{eq:CP}.  Assume that the functions $f_i$ are bounded from below and that Assumption \ref{assump:1} holds. Then, $x(t)$ converges weakly to a weakly Pareto optimal solution of \eqref{eq:MOP}.
    \label{thm:weak_conv}
\end{theorem}

\begin{proof}
    Define the set 
    \begin{align*}
        S \coloneqq \{ z \in \H : f_i(z) \le f_i^{\infty}\, \text{ for all }i = 1,\dots, m\},
    \end{align*}
    where $f_i^{\infty} = \lim_{t\to + \infty} f_i(x(t))$. This limit exists due to Theorem \ref{thm:conv_f_i}. Since $x(t)$ is bounded it posses a weak sequential cluster point $x^{\infty} \in \H$. Hence, there exists a sequence $(x(t_k))_{k \ge 0}$ with $t_k \to + \infty$ and $x(t_k) \rightharpoonup x^{\infty}$ for $k \to + \infty$. Because the objective functions are lower semicontinuous in the weak topology, we get for all $i = 1,\dots, m$
    \begin{align*}
        f_i(x^{\infty}) \le \liminf_{k \to + \infty} f_i(x(t_k)) = \lim_{k \to + \infty} f_i(x(t_k))  = f_i^{\infty}.
    \end{align*}
    Therefore, we can conclude that $x(t)$ converges weakly to an element $x^{\infty} \in S$. Hence, $S$ is nonempty and each weak sequenial cluster point of $x(t)$ belongs to $S$. Let $z \in S$ and define $h_z(t) = \frac{1}{2}\lVert x(t) - z \rVert^2$. The first and second derivative of $h_z(t)$ are given by
    \begin{align*}
        \dot{h}_z(t) = \langle x(t) - z , \dot{x}(t) \rangle \text{ and } \ddot{h}_z(t) = \langle x(t) - z , \ddot{x}(t) \rangle + \lVert \dot{x}(t) \rVert^2,
    \end{align*}
    for almost all $t \in [t_0, +\infty)$. Multiplying $\dot{h}_z(t)$ with $\frac{\alpha}{t}$ and adding it to $\ddot{h}_z(t)$ gives
    \begin{align}
        \ddot{h}_z(t) + \frac{\alpha}{t}\dot{h}_z(t) = \left\langle x(t) - z, \ddot{x}(t) + \frac{\alpha}{t} \dot{x}(t) \right\rangle + \lVert \dot{x}(t) \rVert^2.
        \label{eq:weak_conv_4}
    \end{align}
    Using the equation \eqref{eq:MAVD} together with the weights $\theta(t) \in \Delta^m$ from Lemma \ref{lem:theta_measurable}, we get from \eqref{eq:weak_conv_4} the equation
    \begin{align}
        \ddot{h}_z(t) + \frac{\alpha}{t}\dot{h}_z(t) = \sum_{i=1}^m \theta_i(t)\left\langle z - x(t), \nabla f_i(x(t)) \right\rangle + \lVert \dot{x}(t) \rVert^2.
        \label{eq:weak_conv_4.5}
    \end{align}
    We want to bound the inner products $\left\langle z - x(t), \nabla f_i(x(t)) \right\rangle$. Since $\W_i(t)$ is monotonically decreasing by Proposition \ref{prop:dissipative_const} and converging to $f_i^{\infty}$ by Theorem \ref{thm:conv_f_i}, we get
    \begin{align}
        f_i(x(t)) + \frac{1}{2} \lVert \dot{x}(t) \rVert^2 = \W_i(t) \ge f_i^{\infty},
        \label{eq:weak_conv_5}
    \end{align}
    for all $i = 1,\dots, m$. From $z \in S$ and the convexity of the functions $f_i$, we conclude
    \begin{align}
        f_i^{\infty} \ge f_i(z) \ge f_i(x(t)) + \langle \nabla f_i(x(t)), z - x(t) \rangle.
        \label{eq:weak_conv_6}
    \end{align}
    Together, \eqref{eq:weak_conv_5} and \eqref{eq:weak_conv_6} imply
    \begin{align}
        \langle \nabla f_i(x(t)), z - x(t) \rangle \le \frac{1}{2}\lVert \dot{x}(t) \rVert^2,
        \label{eq:weak_conv_7}
    \end{align}
    for all $i = 1,\dots, m$. Now, we combine \eqref{eq:weak_conv_4.5} and \eqref{eq:weak_conv_7} and multiply with $t$, to conclude
    \begin{align}
        t\ddot{h}_z(t) + \alpha \dot{h}_z(t) \le \frac{3t}{2}\lVert \dot{x}(t) \rVert^2.
        \label{eq:weak_conv_9}
    \end{align}
    Theorem \ref{thm:u_0_x(t)_=_O(t^-2)} states that $(t \mapsto t\lVert\dot{x}(t) \rVert^2) \in L^1([t_0, +\infty))$ for $\alpha > 3$. Then, Lemma \ref{lem:differential_inequality} applied to equation \eqref{eq:weak_conv_9} guarantees that $h_z(t)$ converges and by Opial's lemma (Lemma \ref{lem:opial}) we conclude that $x(t)$ converges weakly to an element in $S$. By Theorem \ref{thm:conv_W}, we know that every weak accumulation point of $x(t)$ is weakly Pareto optimal.
\end{proof}

\section{Numerical experiments}
\label{sec:numerical_experiments}
In this section, we conduct numerical experiments to verify the convergence rates we prove in the previous section. In particular, we show that the convergence of $u_0(x(t))$ with rate $\mathcal{O}(t^{-2})$ as stated in Theorem \ref{thm:u_0_x(t)_=_O(t^-2)} holds. Since we cannot calculate analytical solutions to \eqref{eq:CP} for a general multiobjective optimization problem in closed form, we compute the approximation to a solution $x$ using a discretization. We do not discuss the quality of the discretization we use. For all experiments we use initial time $t_0 = 1$, set a fixed initial state $x(t_0) = x_0$ and use initial velocity $\dot{x}(t_0) = 0$. We use equidistant time steps $t_k = t_0 + kh$, with $h = 1\mathrm{e}{-3}$. We use the scheme $x(t_k) \approx x^k$, $\dot{x}(t_k) \approx \frac{x^{k+1} - x^k}{h}$ and $\ddot{x}(t_k) \approx \frac{x^{k+1} - 2x^k + x^{k-1}}{h^2}$ to compute the discretization $(x^k)_{k \ge 0}$ of the trajectory $x$ for $100\,\,000$ time steps. We look at two example with instances of the multiobjective optimization problem \eqref{eq:MOP}. Both problem instances use two convex and smooth objective functions $f_i: \R^2 \to \R$ for $i=1,2$. In Subsection \ref{subsec:quad_prob} we look at a quadratic multiobjective optimization problem and in Subsection \ref{subsec:non_quad_prob} we consider a convex optimization problem with objective functions that are not strongly convex. For both examples we plot approximations of the solution $x$ and plot the function $u_0(x(t))$ to show that the inequality $u_0(x(t)) \le \frac{t_0^2 u_0(x_0) + 2(\alpha - 1)R}{t^2}$ holds for $t \ge t_0$. To compute $u_0(x(t))$ we have to solve the optimization problem $u_0(x^k) = \sup_{z \in \H} \min_{i = 1,\dots,m} f_i(x^k) - f_i(z)$ for every of the $100\,\,000$ iterations with adequate accuracy. Therefore, we restrict ourselves to problems where the Pareto set of \eqref{eq:MOP} can be explicitly computed. For these problems $u_0$ can be evaluated more efficiently using Lemma \ref{lem:sup_inf_u_0}.

\begin{figure}
    \begin{center}        
        \begin{subfigure}[t]{.24\textwidth}
            \centering
            \includegraphics[width=\linewidth]{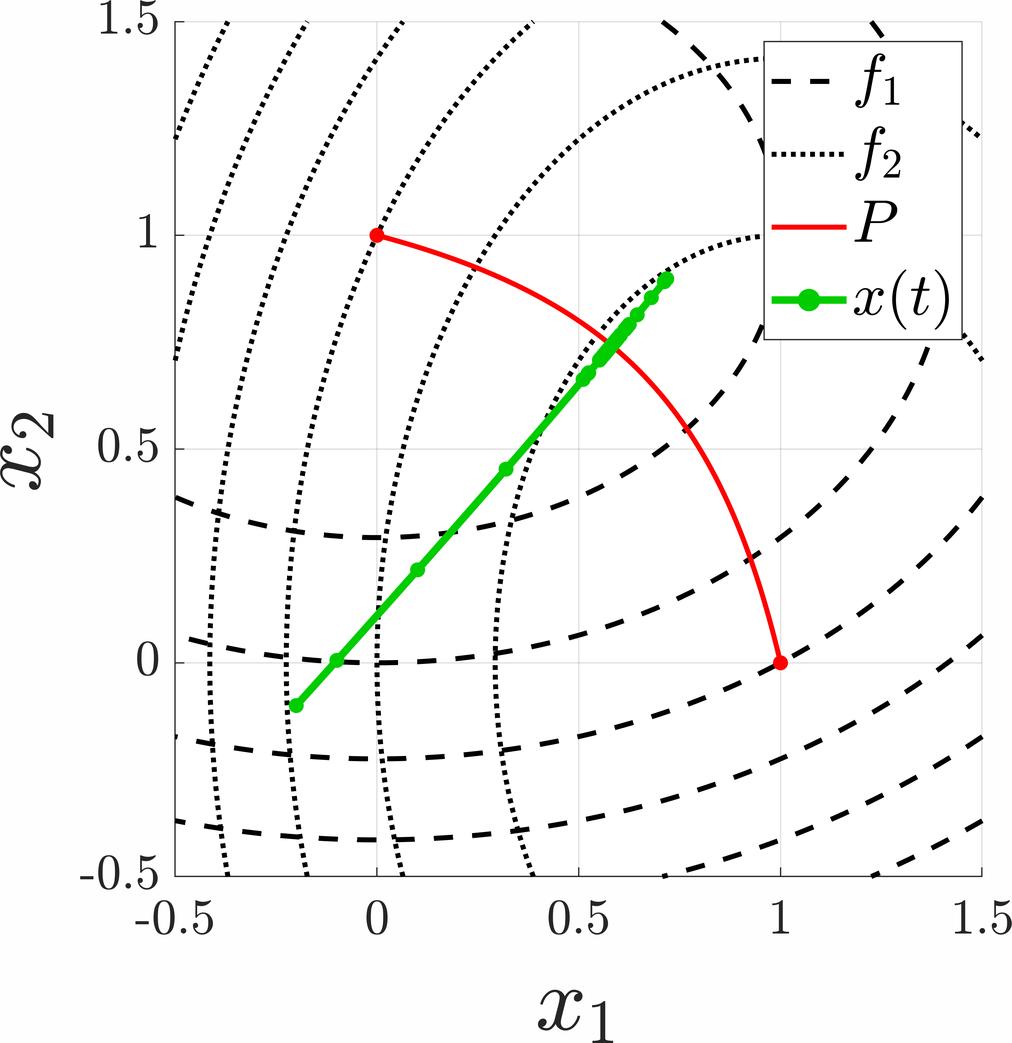}
            \caption{$\alpha=3$}
            \label{subfig:all_quad_a}
        \end{subfigure}
        \begin{subfigure}[t]{.24\textwidth}
            \centering
            \includegraphics[width=\linewidth]{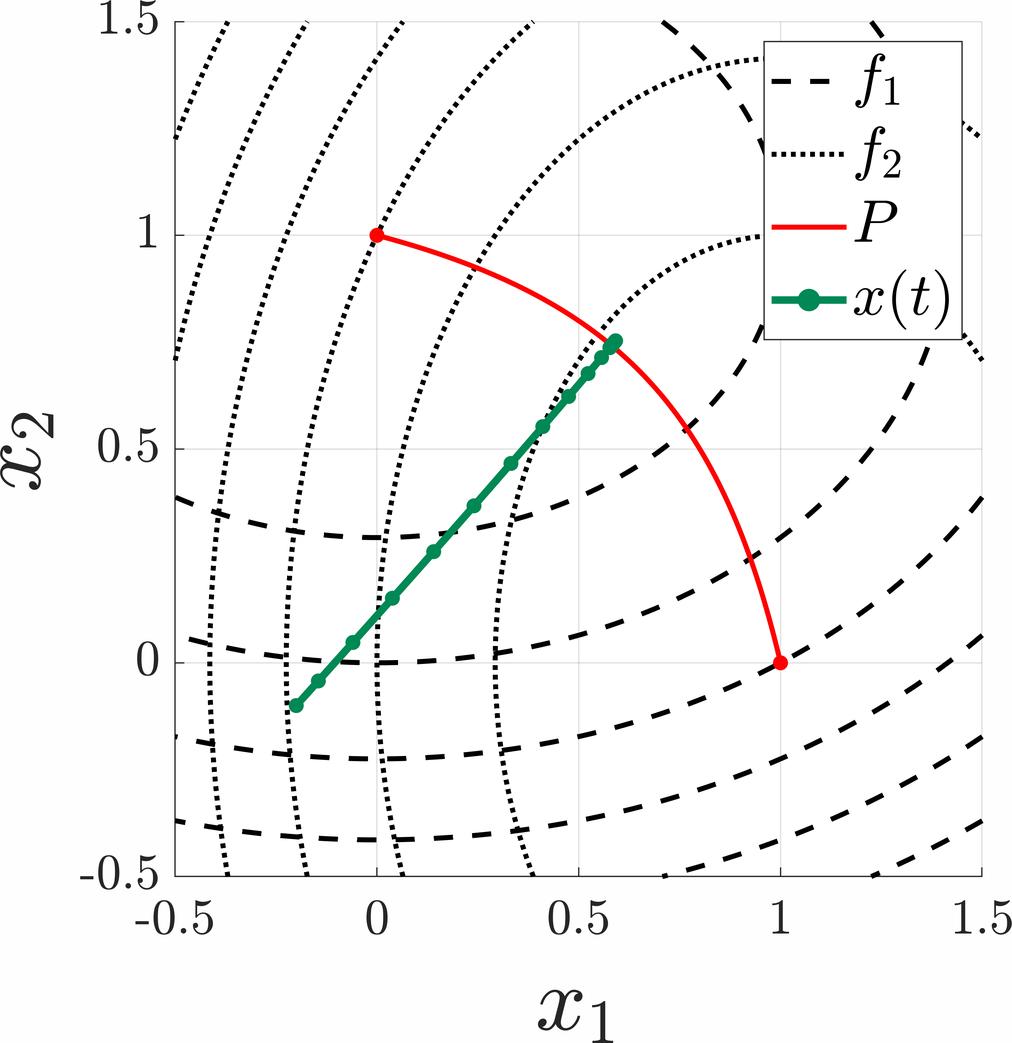}
            \caption{$\alpha=10$}
        \end{subfigure}    
          \begin{subfigure}[t]{.24\textwidth}
            \centering
            \includegraphics[width=\linewidth]{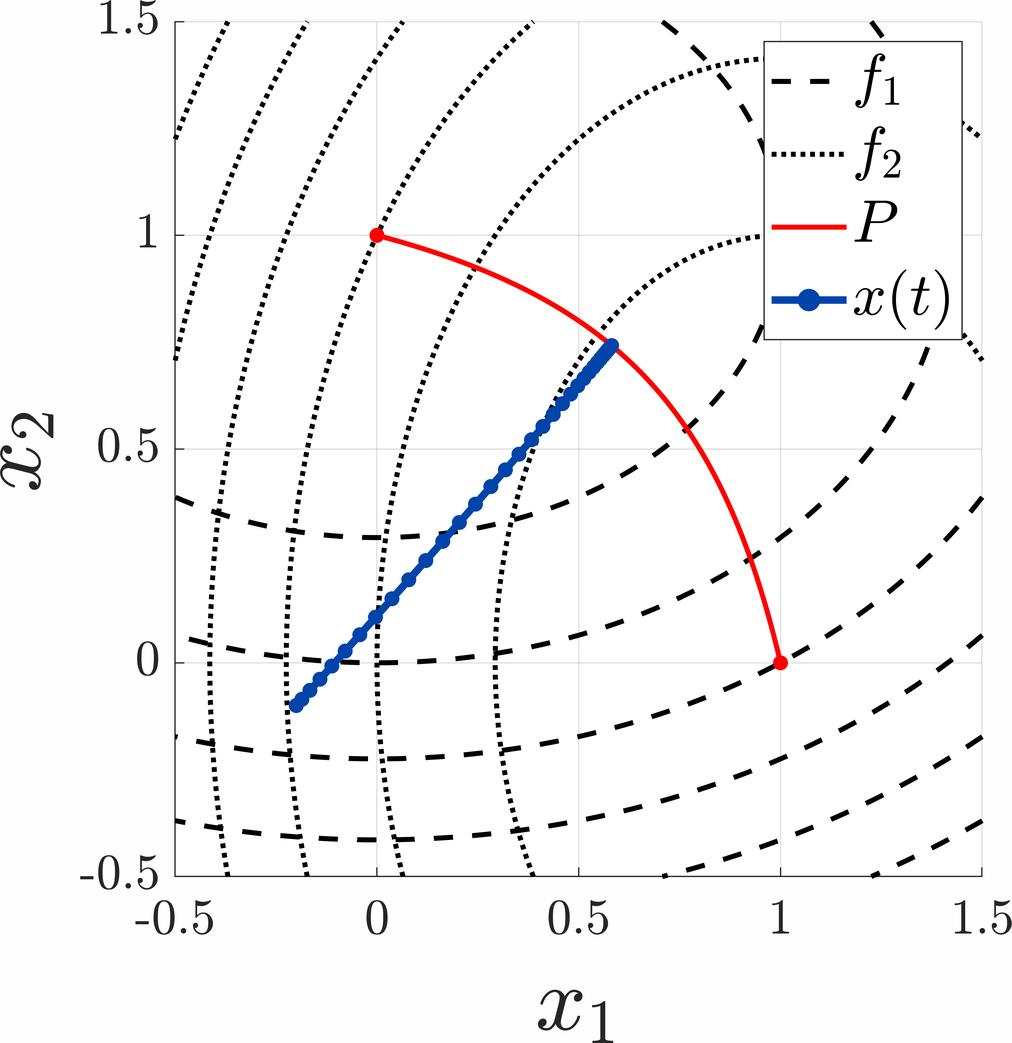}
            \caption{$\alpha=50$}
          \end{subfigure}
          \begin{subfigure}[t]{.24\textwidth}
            \centering
            \includegraphics[width=\linewidth]{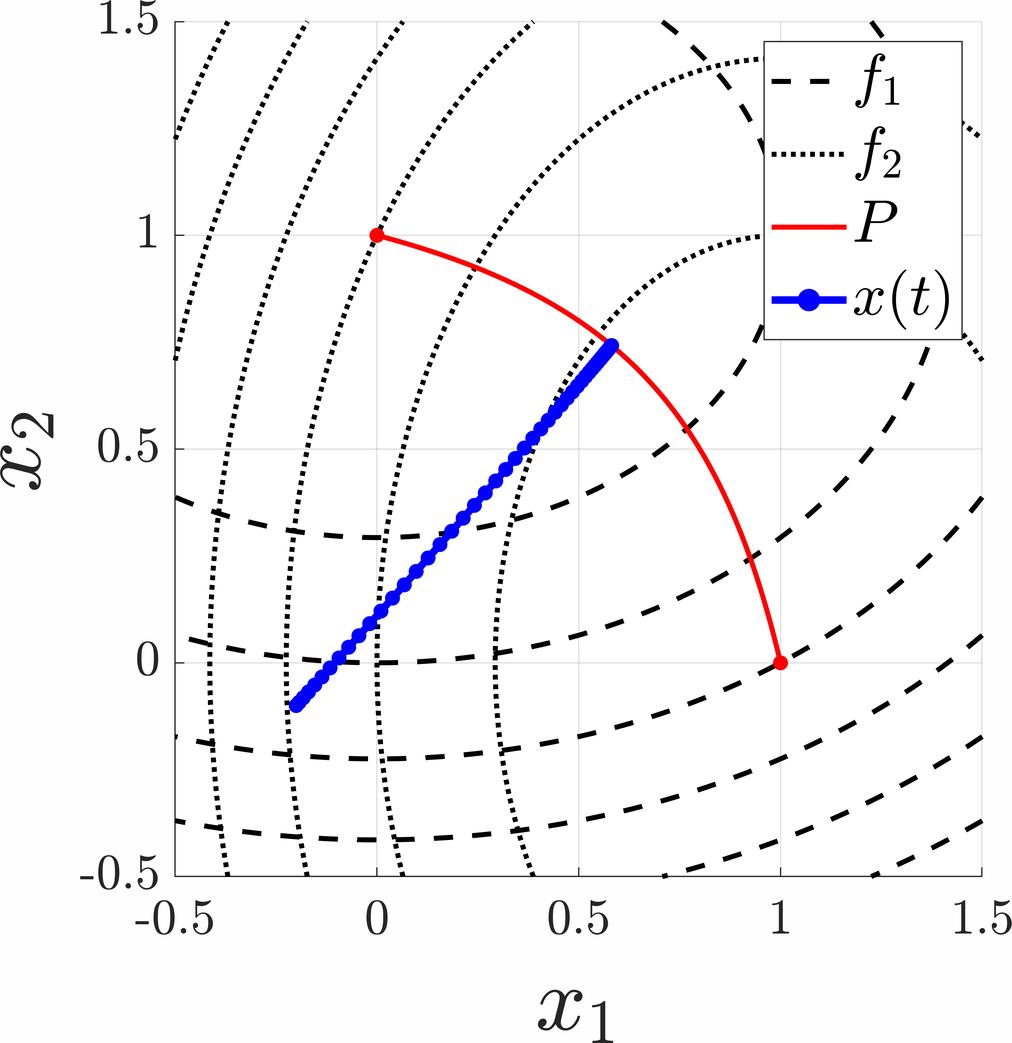}
            \caption{$\alpha=100$}
            \label{subfig:all_quad_d}
          \end{subfigure}
    \medskip
          \begin{subfigure}[t]{.24\textwidth}
            \centering
            \includegraphics[width=\linewidth]{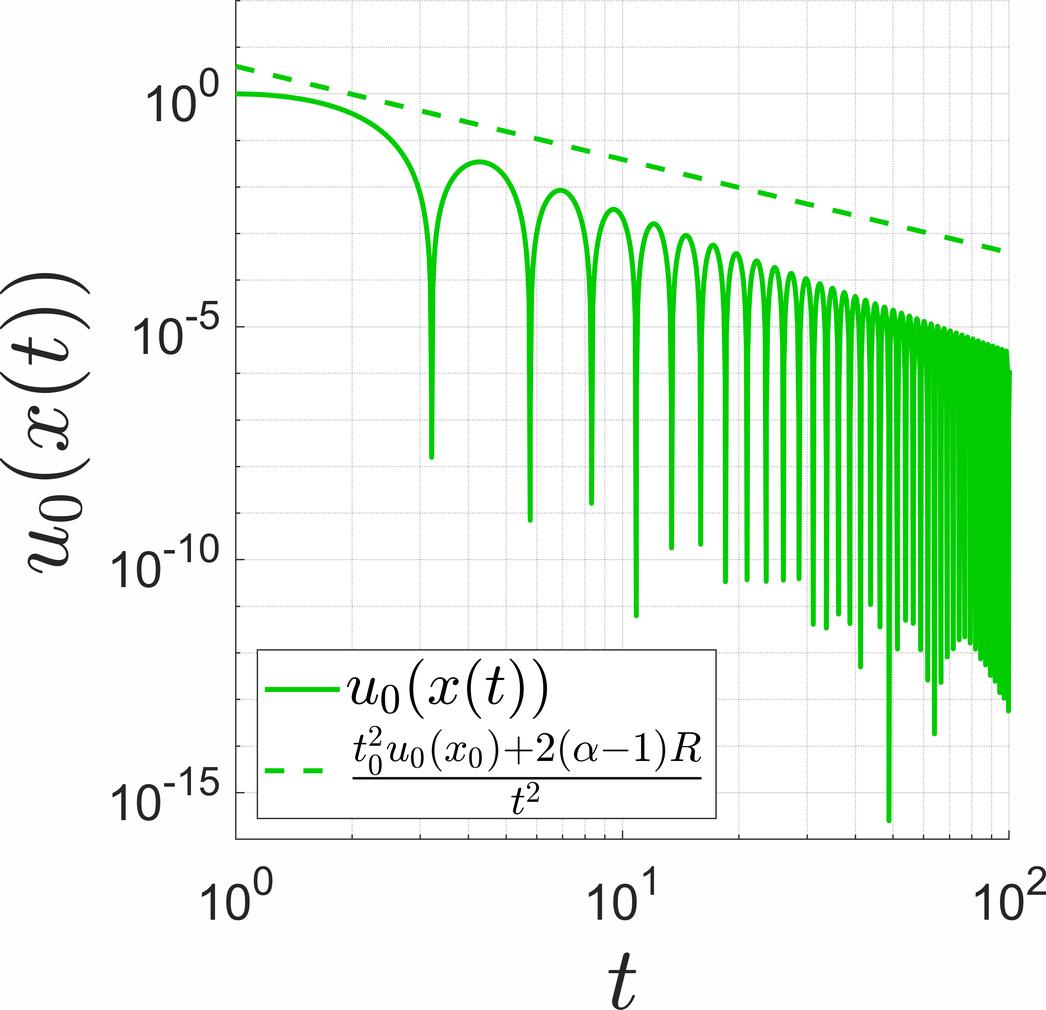}
            \caption{$\alpha=3$}
            \label{subfig:all_quad_e}
          \end{subfigure}
          \begin{subfigure}[t]{.24\textwidth}
            \centering
            \includegraphics[width=\linewidth]{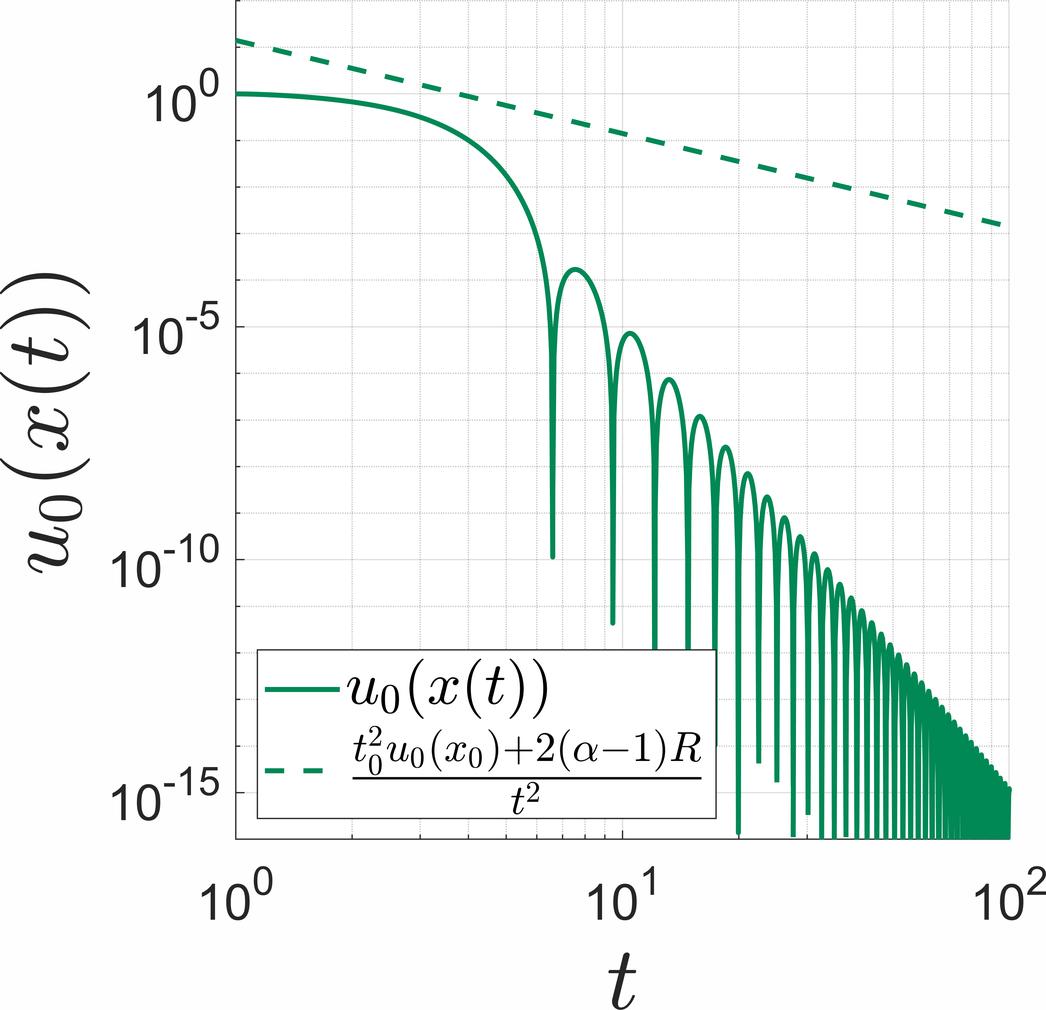}
            \caption{$\alpha=10$}
          \end{subfigure}    
          \begin{subfigure}[t]{.24\textwidth}
            \centering
            \includegraphics[width=\linewidth]{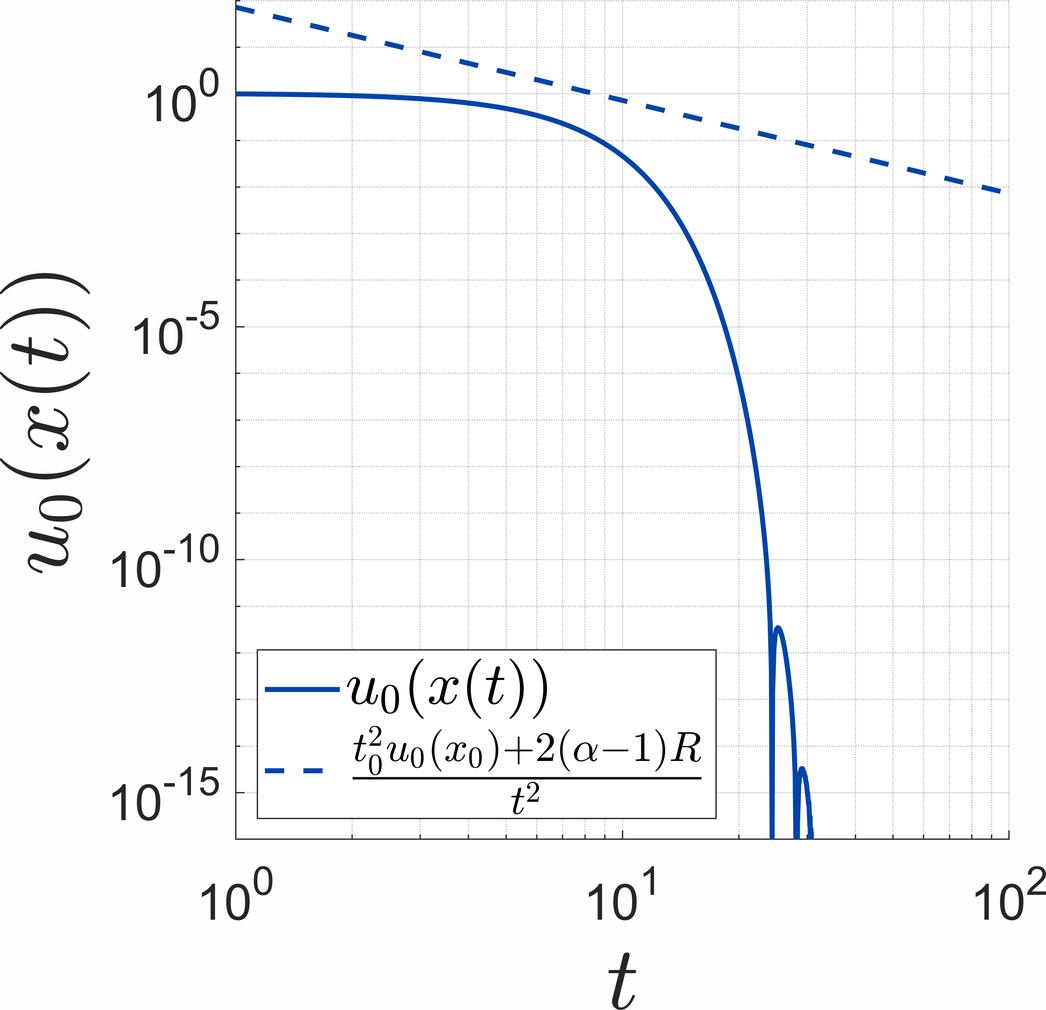}
            \caption{$\alpha=50$}
          \end{subfigure}
          \begin{subfigure}[t]{.24\textwidth}
            \centering
            \includegraphics[width=\linewidth]{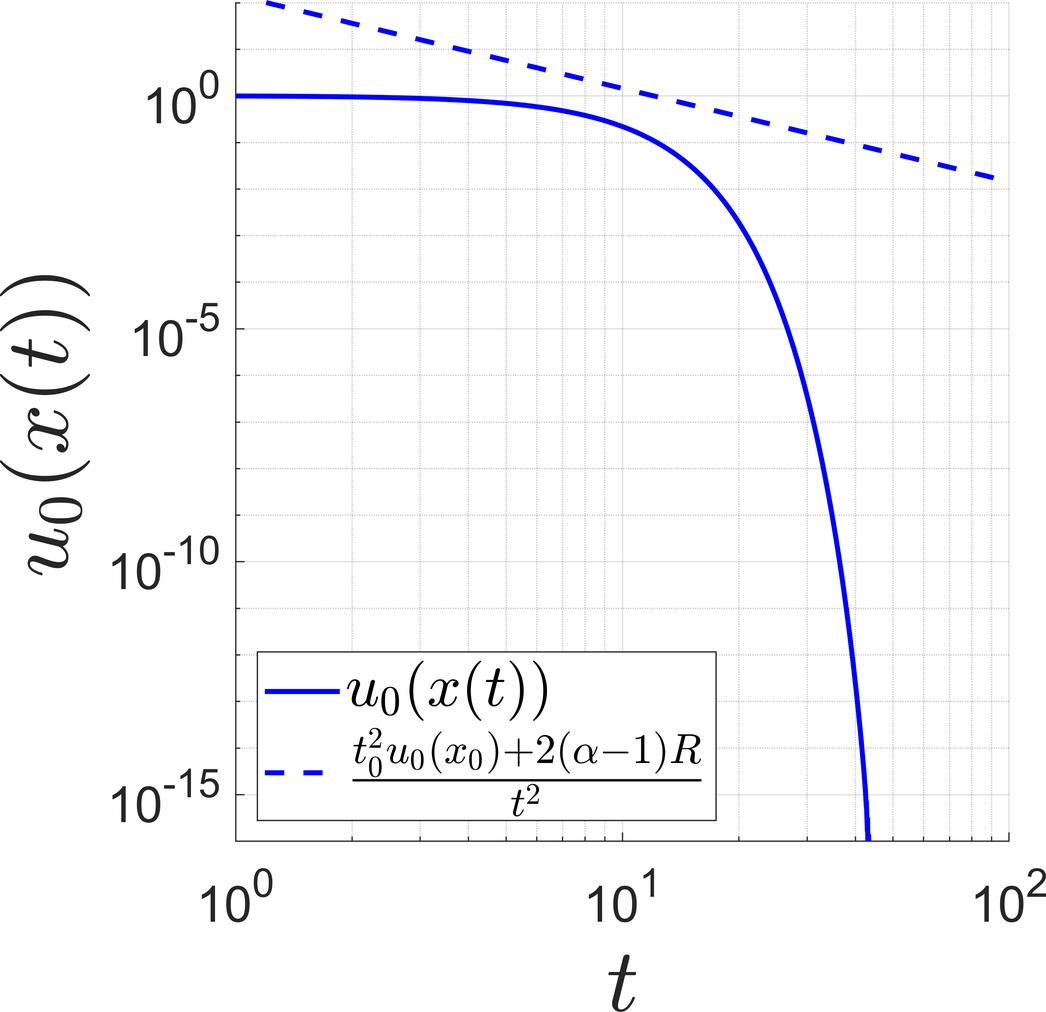}
            \caption{$\alpha=100$}
            \label{subfig:all_quad_h}
          \end{subfigure}
          \caption{Trajectories $x$ and inequalities $u_0(x(t)) \le \frac{t_0^2 u_0(x_0) + 2(\alpha - 1)R}{t^2}$ for different values of $\alpha \in \{3, 10, 50, 100 \}$.}
          \label{fig:all_quad}
    \end{center}
\end{figure}

\subsection{A quadratic multiobjective optimization problem}
\label{subsec:quad_prob}
We begin with an instance of \eqref{eq:MOP} with two quadratic objective functions
\begin{align*}
    f_i: \R^2 \to \R, \quad x \mapsto \frac{1}{2}(x-x^i)^{\top} Q_i (x-x^i),
\end{align*}
for $i = 1,2$, given matrices and vectors
\begin{align*}
    Q_1 = \left( \begin{array}{cc}
        2 & 0 \\
        0 & 1
    \end{array} \right),\quad
    Q_2 = \left( \begin{array}{cc}
        1 & 0 \\
        0 & 2
    \end{array} \right),\quad
    x^1 = \left( \begin{array}{c}
        1  \\
        0
    \end{array} \right),\quad
    x^2 = \left( \begin{array}{c}
        0  \\
        1
    \end{array} \right).
\end{align*}
For this problem the Pareto set is
\begin{align*}
    P = \left\lbrace x \in \R^2 : x = \left( \begin{array}{c}
        2\lambda/(1 + \lambda)  \\
        2(1-\lambda)/(2 - \lambda)
    \end{array} \right), \text{ for } \lambda \in [0,1] \right\rbrace.
\end{align*}
In our first experiment, we use the initial value $x_0 = (-.2, -.1)^{\top}$. We compute an approximation of a solution to \eqref{eq:CP} for different values of $\alpha \in\{3, 10, 50, 100\}$ as described in the introduction of Section \ref{sec:numerical_experiments}. The results can be seen in Figure \ref{fig:all_quad}. Subfigures \ref{subfig:all_quad_a} - \ref{subfig:all_quad_d} contain plots of the trajectories $x$ for different values of $\alpha$. In the plots of the trajectories we added a circle every $500$ iterations to visualize the velocities. In Subfigures \ref{subfig:all_quad_e} - \ref{subfig:all_quad_h} the values of $u_0(x(t))$ and the bounds $\frac{t_0^2 u_0(x_0) + 2(\alpha - 1)R}{t^2}$ for different values of $\alpha$ are shown. The inequality $u_0(x(t)) \le \frac{t_0^2 u_0(x_0) + 2(\alpha - 1)R}{t^2}$ holds for each value of $\alpha$. For the smallest value of $\alpha=3$ we see a large number of oscillations in the trajectory and in the values of $u_0(x(t))$, respectively. This behavior is typical for systems with asymptotic vanishing damping. For larger values of $\alpha$ we observe fewer oscillations and see improved convergence rates, with slower movement in the beginning due to the high friction. These phenomena are consistent with the observations made in the singleobjective setting.

\begin{figure}
    \begin{center}
        \begin{subfigure}[t]{.24\textwidth}
            \centering
            \includegraphics[width=\linewidth]{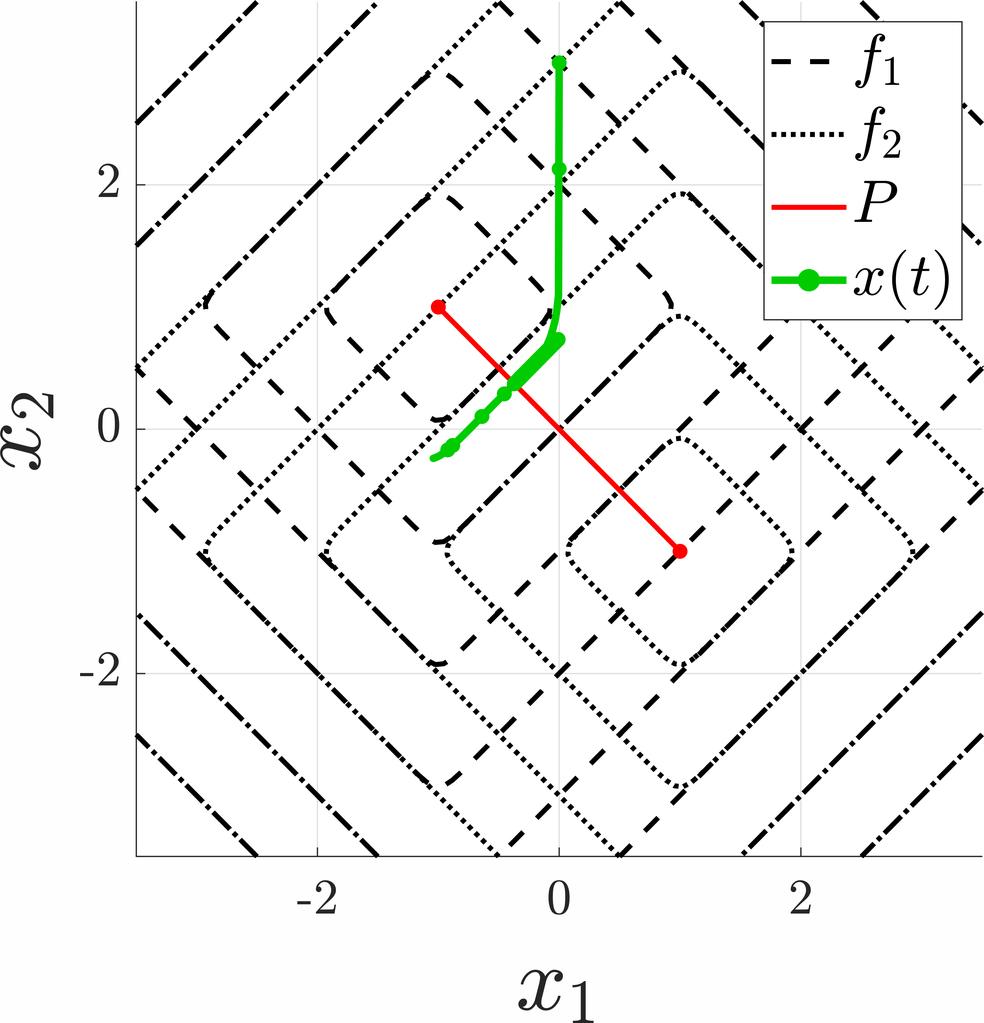}
            \caption{$\alpha=3$}
            \label{subfig:all_nonquad_a}
        \end{subfigure}
        \begin{subfigure}[t]{.24\textwidth}
            \centering
            \includegraphics[width=\linewidth]{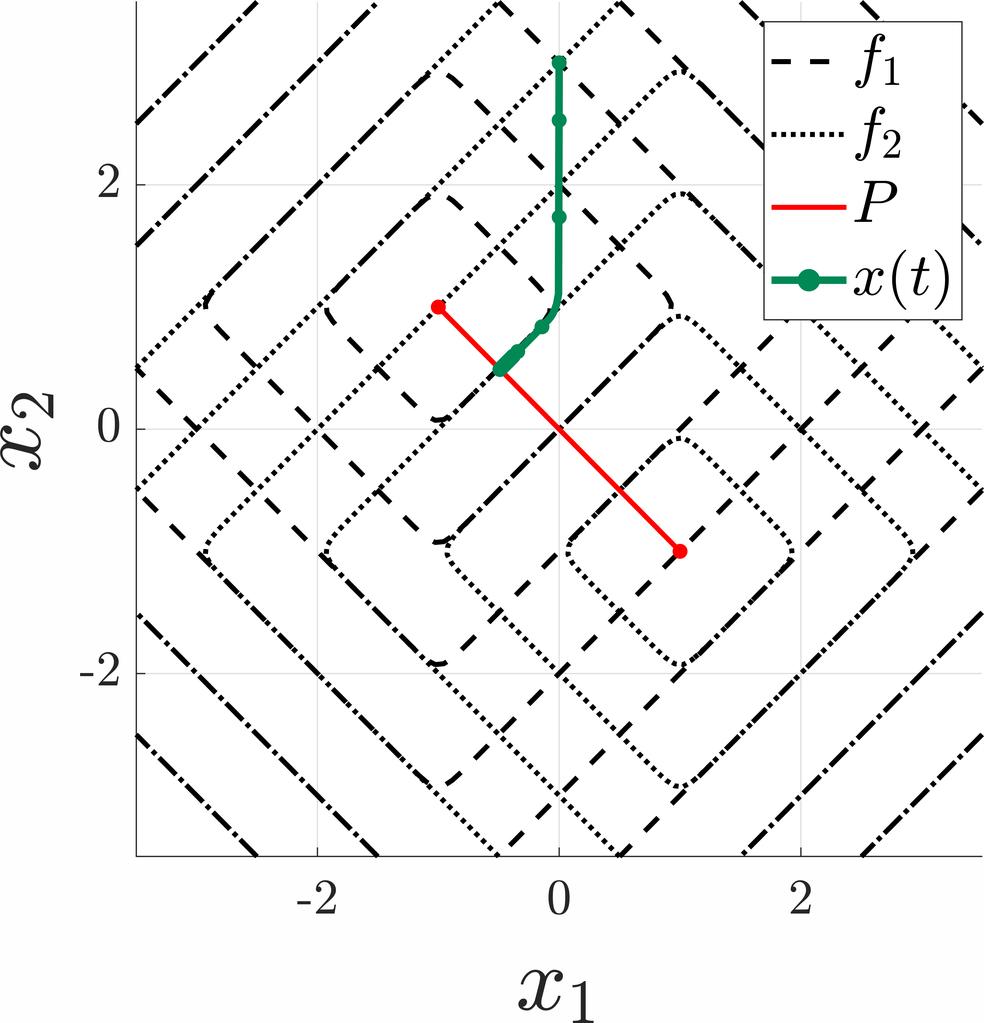}
            \caption{$\alpha=10$}
        \end{subfigure}    
          \begin{subfigure}[t]{.24\textwidth}
            \centering
            \includegraphics[width=\linewidth]{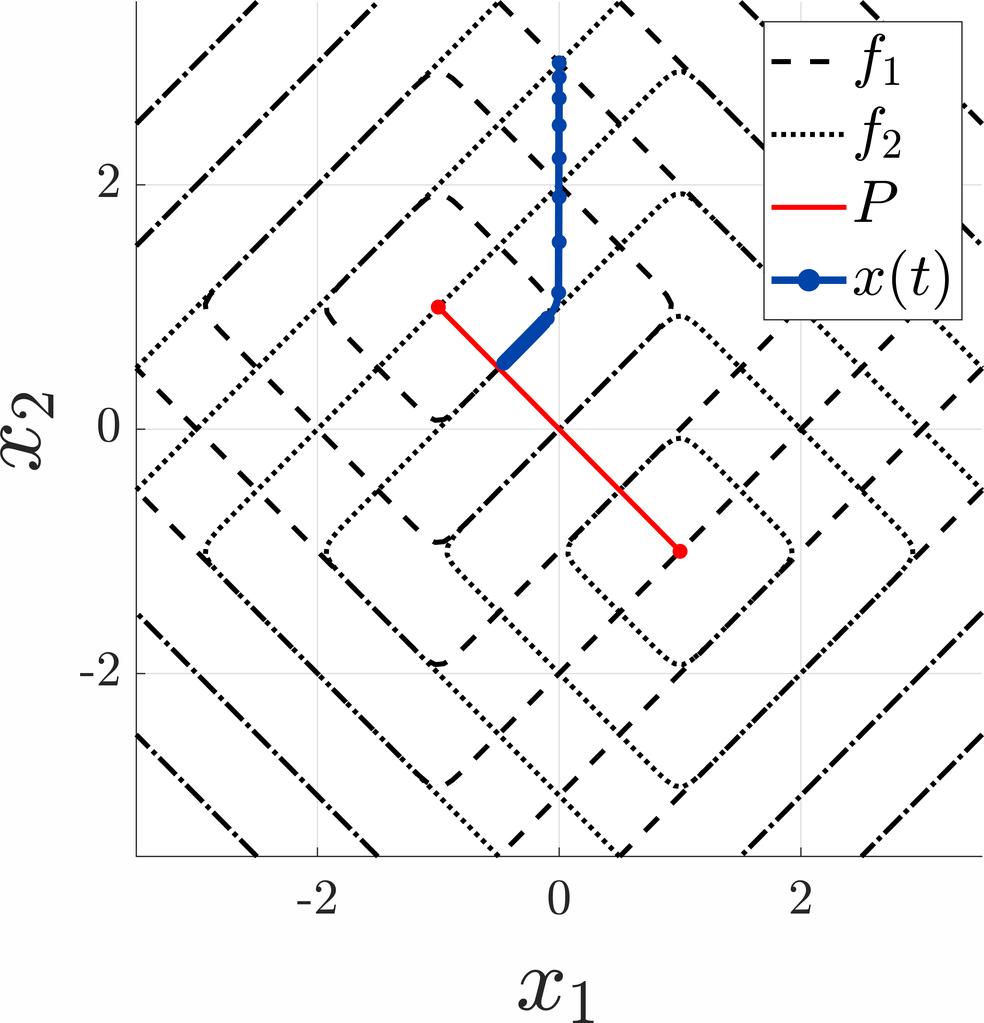}
            \caption{$\alpha=50$}
          \end{subfigure}
          \begin{subfigure}[t]{.24\textwidth}
            \centering
            \includegraphics[width=\linewidth]{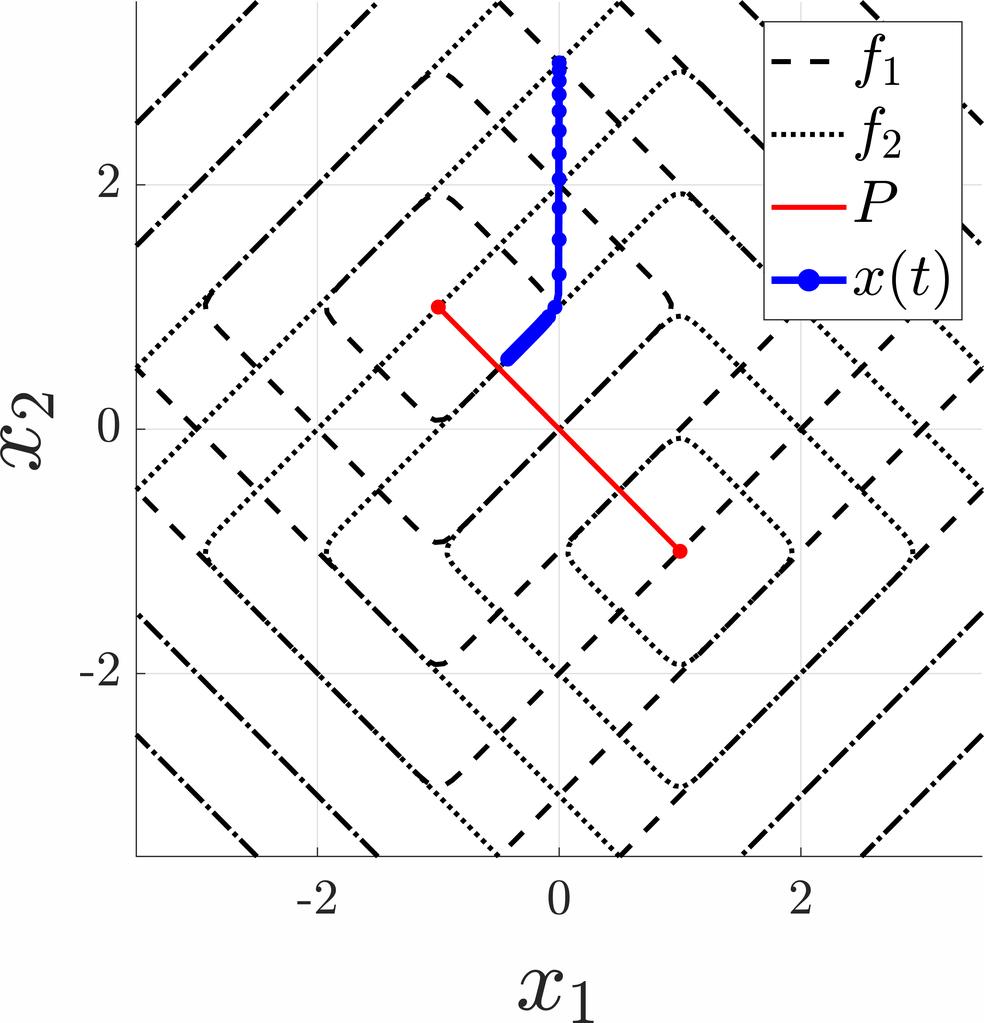}
            \caption{$\alpha=100$}
            \label{subfig:all_nonquad_d}
          \end{subfigure}
    \medskip
          \begin{subfigure}[t]{.24\textwidth}
            \centering
            \includegraphics[width=\linewidth]{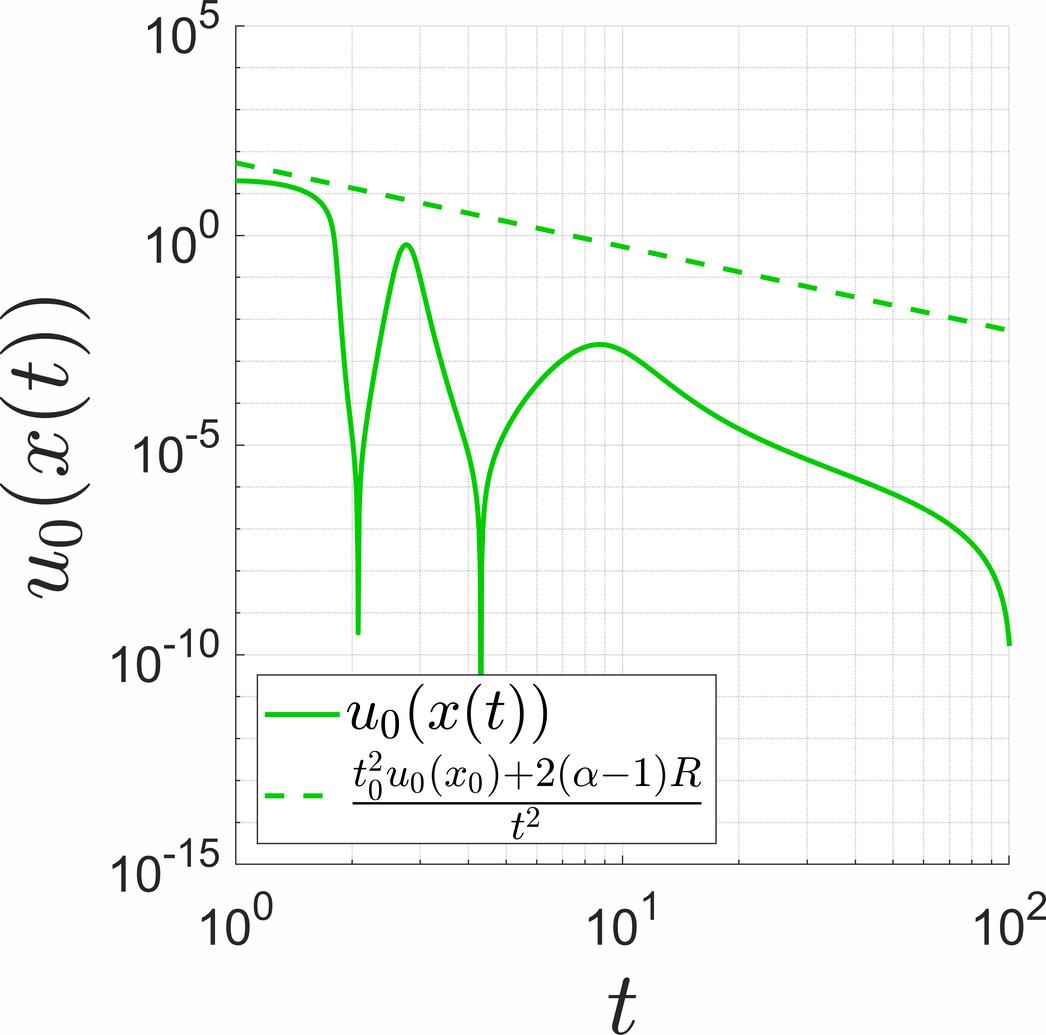}
            \caption{$\alpha=3$}
            \label{subfig:all_nonquad_e}
          \end{subfigure}
          \begin{subfigure}[t]{.24\textwidth}
            \centering
            \includegraphics[width=\linewidth]{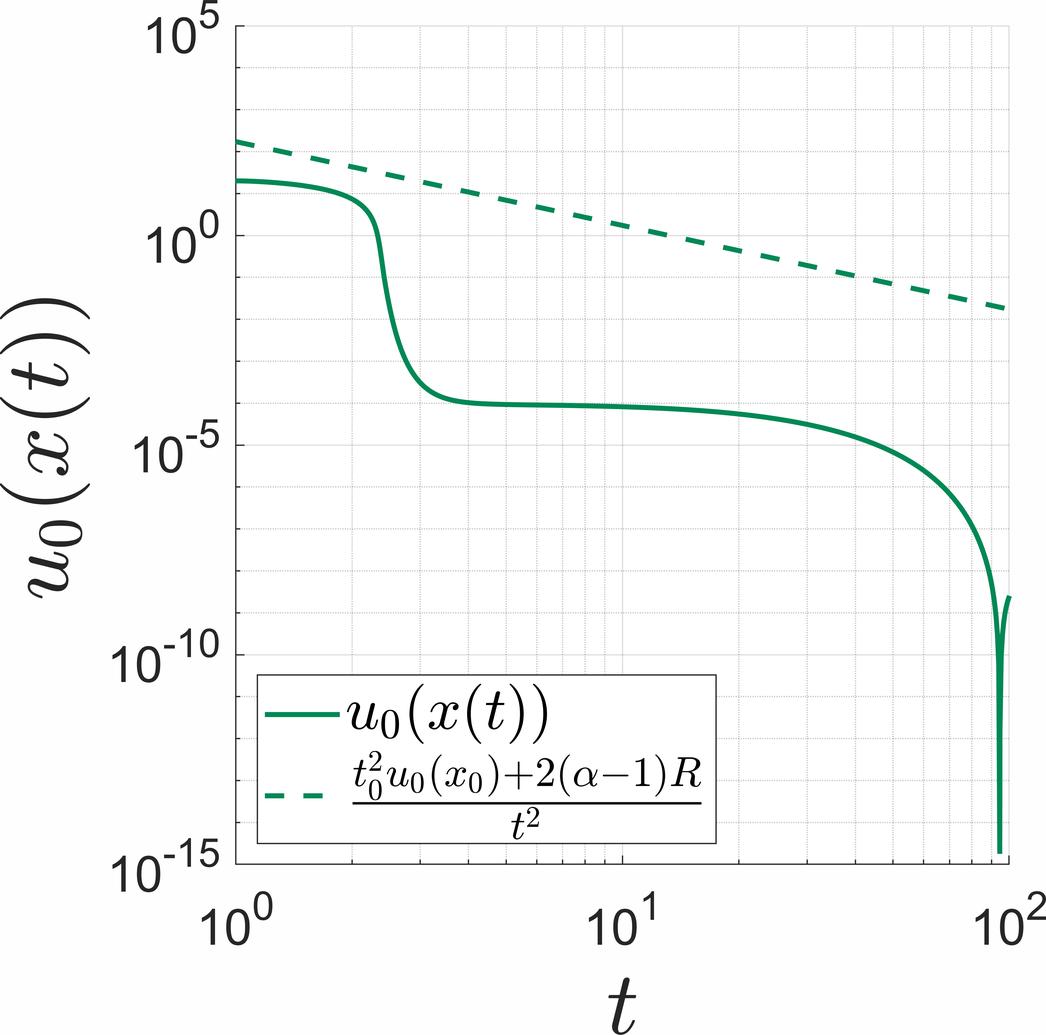}
            \caption{$\alpha=10$}
          \end{subfigure}    
          \begin{subfigure}[t]{.24\textwidth}
            \centering
            \includegraphics[width=\linewidth]{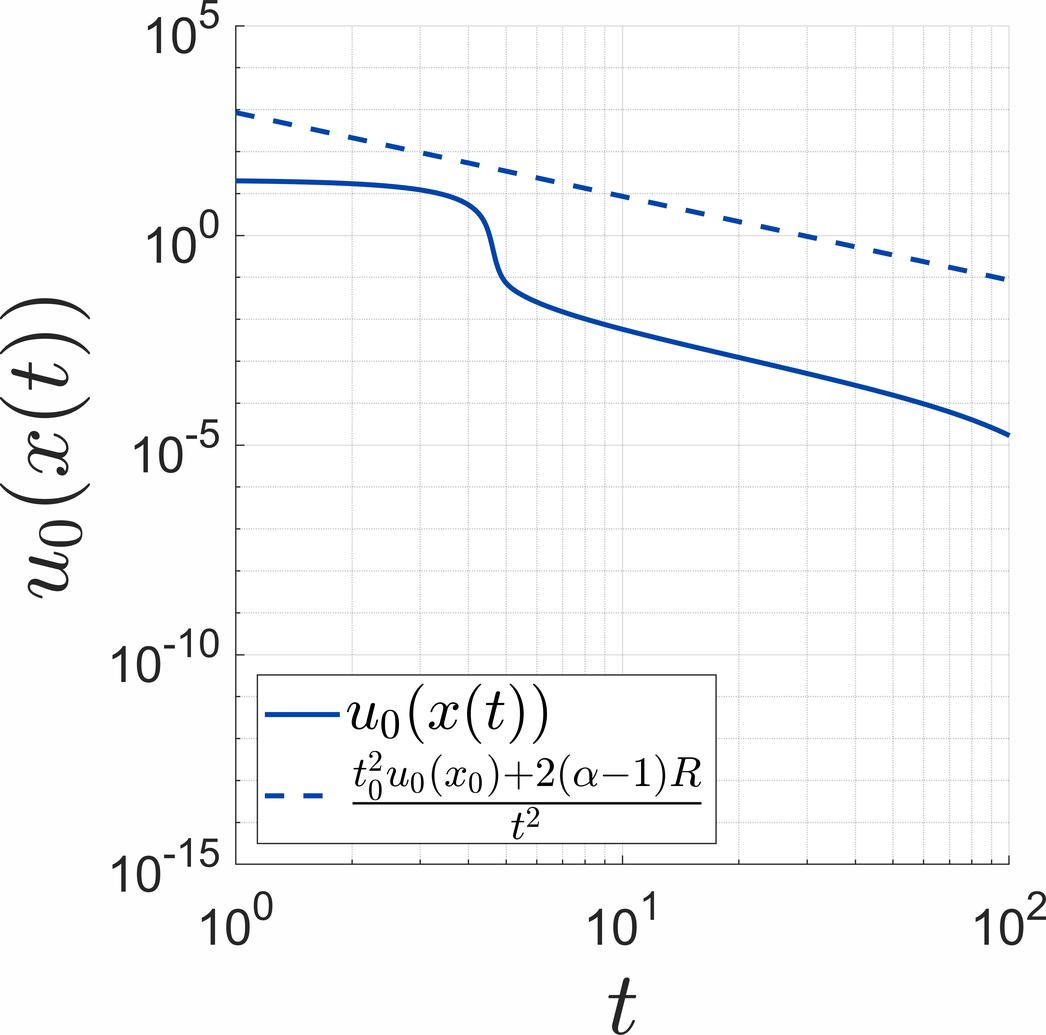}
            \caption{$\alpha=50$}
          \end{subfigure}
          \begin{subfigure}[t]{.24\textwidth}
            \centering
            \includegraphics[width=\linewidth]{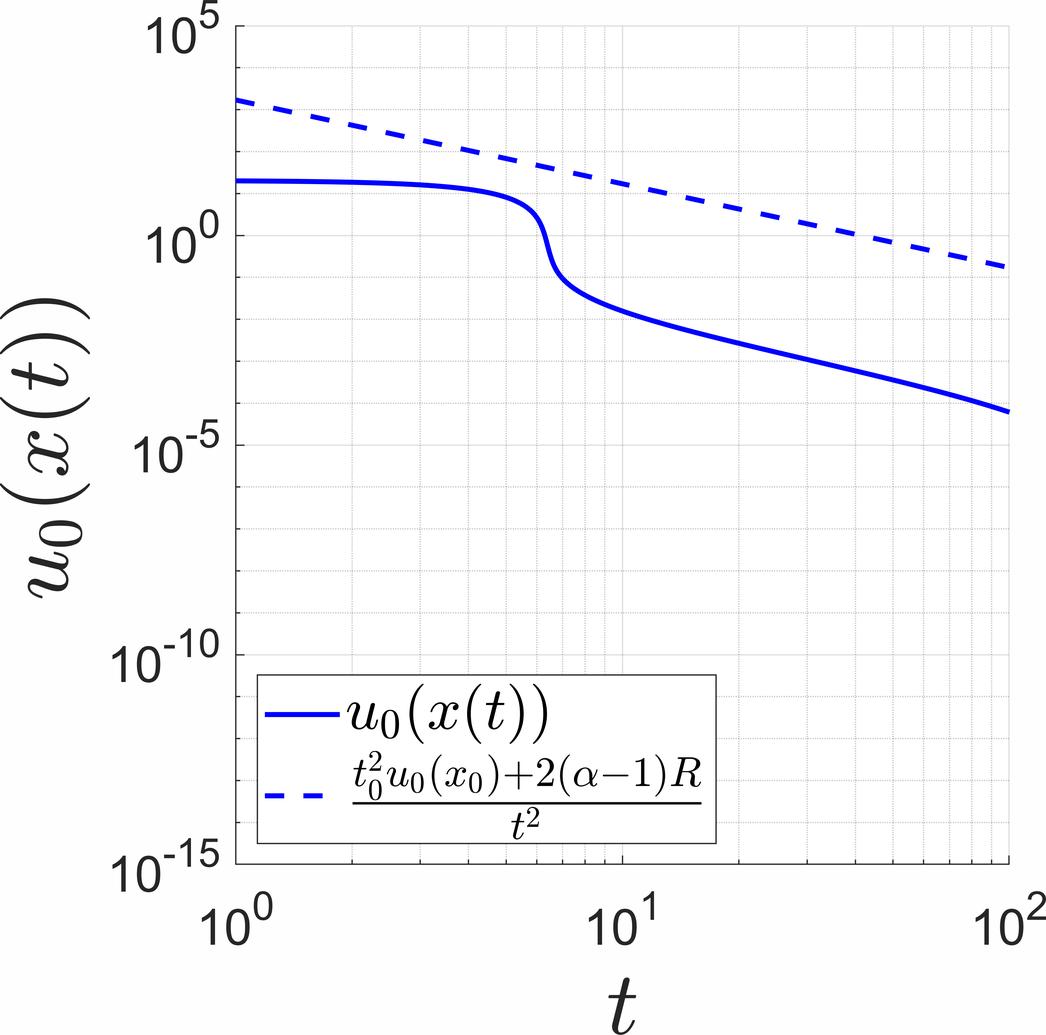}
            \caption{$\alpha=100$}
            \label{subfig:all_nonquad_h}
          \end{subfigure}
          \caption{Trajectories $x$ and inequalities $u_0(x(t)) \le \frac{t_0^2 u_0(x_0) + 2(\alpha - 1)R}{t^2}$ for different values of $\alpha \in \{3, 10, 50, 100 \}$.}
           \label{fig:all_nonquad}
    \end{center}
\end{figure}
\subsection{A nonquadratic multiobjective optimization problem}
\label{subsec:non_quad_prob}
In our second example, we consider the problem \eqref{eq:MOP} with two objective functions
\begin{align}
\label{eq:obj_func_nonquad}
    f_i: \R^2 \to \R, \quad x \mapsto \log\left( \sum_{j=1}^p \exp\left(\left(a_j^{(i)}\right)^\top x - b_j^{(i)}\right)\right),
\end{align}
for $i = 1,2$, $p = 4$ and given matrices and vectors
\begin{align*}
    A^{(1)} = 
    \left( \begin{array}{c}
        \left(a_1^{(1)}\right)^\top \\
        \vdots\\
        \left(a_4^{(1)}\right)^\top
    \end{array} \right) 
    =
    \left( \begin{array}{rr}
        10 & 10 \\
        10 & -10 \\
        -10 & -10 \\
        -10 & 10
    \end{array} \right),\quad
    b^{(1)} = \left( \begin{array}{r}
        0 \\
        -20 \\
        0 \\
        20
    \end{array} \right),\\
    A^{(2)} =
    \left( \begin{array}{c}
        \left(a_1^{(2)}\right)^\top \\
        \vdots\\
        \left(a_4^{(2)}\right)^\top
    \end{array} \right) 
    =
    \left( \begin{array}{rr}
        10 & 10 \\
        10 & -10 \\
        -10 & -10 \\
        -10 & 10
    \end{array} \right),\quad
    b^{(2)} = \left( \begin{array}{r}
        0 \\
        20 \\
        0 \\
        -20
    \end{array} \right).
\end{align*}
The objective functions given by \eqref{eq:obj_func_nonquad} are convex but not strongly convex. Taking advantage of the symmetry in the objective functions $f_i$, the Pareto set $P$ can be explicitly computed and is
\begin{align*}
    P = \left\lbrace x \in \R^2 : x = \left( \begin{array}{r}
        -1 + 2\lambda \\
        1 - 2\lambda
    \end{array} \right), \text{ for } \lambda \in [0,1] \right\rbrace.
\end{align*}
We choose the initial value $x_0 = (0,3)^{\top}$ and compute an approximate solution to \eqref{eq:CP} as described in the beginning of Section \ref{sec:numerical_experiments}. Analogous to the last example, we present the results of the computations in Figure \ref{fig:all_nonquad}. Again, Subfigures \ref{subfig:all_nonquad_a} - \ref{subfig:all_nonquad_d} contain plots of the trajectories and Subfigures \ref{subfig:all_nonquad_e} - \ref{subfig:all_nonquad_h} contain the values of the merit function $u_0(x(t))$. We observe results similar to the example in Subsection \ref{subsec:quad_prob}. Since the objective functions given in \eqref{eq:obj_func_nonquad} are not strongly convex, we experience slower convergence especially in the beginning, where the gradients along the trajectories remain almost constant. Once more, we see for small values of $\alpha$ oscillations in the trajectory $x$ and the merit function values $u_0(x(t))$ introduced by the inertia in the system \eqref{eq:CP}. Larger values of $\alpha$ correspond to higher friction in the beginning and we therefore experience slower convergence for the time interval we consider. Oscillations can only be seen for $\alpha=3$ and close to the end for $\alpha=10$. The slower convergence in this example is expected due to the lack of strong convexity.

\section{Conclusion}
\label{sec:conclusion}
We introduce the system \eqref{eq:MAVD} and discuss its main properties. To the best of our knowledge, this system is the first inertial gradient-like system with asymptotic vanishing damping for multiobjective optimization problems, expanding the ideas laid out in \cite{Attouch2015}. We prove existence of global solutions in finite dimensions for arbitrary initial conditions. We discuss the asymptotic behaviour of solutions to \eqref{eq:MAVD} and show that the function values decrease along the trajectory with rate $u_0(x(t)) = \mathcal{O}(t^{-2})$ for $\alpha \ge 3$ using a Lyapunov type analysis. Further, we show that bounded solutions converge weakly to weakly Pareto optimal points given $\alpha > 3$. These statements are consistent with the result obtained for singleobjective optimization. We verify our results on two test problems and show that the given bounds on the decay of the function values are satisfied.

For future work it would be interesting to further adapt the system \eqref{eq:MAVD}. Possible research directions involve Tikhonov regularization \cite{Attouch2022_2}, Hessian-driven damping \cite{Attouch2023} and the treatment of linear constraints by the means of Augmented Lagrangian type systems \cite{Bot2022}. From an algorithmic point of view it would be interesting to further analyze numerical schemes for fast multiobjective optimization based on \eqref{eq:MAVD} and also to define improved algorithms based on Tikhonov regularization, Hessian-driven damping and for problems with linear constraints.

\nocite{Fliege2000}
\nocite{Mukai1980}
\nocite{attouch2022}
\bibliographystyle{siamplain}
{\footnotesize
\bibliography{references}
}
\section*{Appendix}
\appendix


\section{Proof of Proposition \ref{prop:set_valued_G}}
\label{sec:appendix_set_valued_G}
Here, we recall the most important definitions on set-valued maps we need to prove Proposition \ref{prop:set_valued_G}. The notation is aligned with \cite{Aubin2012}. Let $\mathcal{X}, \mathcal{Y}$ be real Hilbert spaces and let $G : \mathcal{X} \rightrightarrows \mathcal{Y}$ be a set-valued map.
\begin{mydef}
    We say $G$ is upper semicontinuous (u.s.c.) at $x_0 \in \mathcal{X}$ if for any open set $N \subset \mathcal{Y}$ containing $G(x_0)$ there exists a neighborhood $M \subset \mathcal{X}$ of $x_0$ such that $G(M) \subset N$.

    We say that $G$ is u.s.c. if it is so at every $x_0 \in \mathcal{X}$.
\end{mydef}

\begin{mydef}
    We say $G$ is u.s.c. at $x_0$ in the $\varepsilon$ sense if, given $\varepsilon > 0$, there exists $\delta > 0$ such that $G(B_{\delta}(x^0)) \subset G(x^0) + B_{\varepsilon}(0)$.

    We say that $G$ is u.s.c. in the $\varepsilon$ sense if it is so at every $x_0 \in \mathcal{X}$.
\end{mydef}

\begin{myprop}
    Let $G$ be a set valued map. The following statements hold.
    \begin{enumerate}[(i)]
        \item If $G$ is u.s.c. it is also u.s.c. in the $\varepsilon$ sense.
        \item If $G$ is u.s.c. in the $\varepsilon$ sense and takes compact values $G(x) \subset \mathcal{Y}$ for all $x \in \mathcal{X}$, then it is u.s.c. as well.
    \end{enumerate}
    \label{prop:usc_varepsilon_sense}
\end{myprop}

\begin{mydef}
    We say that a map $\phi: \mathcal{X} \to \mathcal{Y}$ is locally compact if for each point $x_0 \in \mathcal{X}$ there exists a neighborhood which is mapped into a compact subset of $\mathcal{Y}$.
\end{mydef}

In the proof of Proposition \ref{prop:set_valued_G} we use the following auxiliary lemma. Lemma \ref{lem:lin_opt_pert} states that the set-valued map $(u,v) \mapsto \argmin_{g \in C(u)} \langle g, -v \rangle$ is u.s.c..

\begin{mylemma}
Let $( \overline{u}, \overline{v}) \in \H \times \H$ be fixed. Then, for all $\varepsilon > 0$ there exists $\delta > 0$ such that for all $(u,v) \in \H\times\H$ with $\lVert (u,v) - (\overline{u}, \overline{v}) \rVert_{\H \times \H} < \delta$ and for all $g \in \argmin_{g \in C(u)} \langle g, -v \rangle$ there exists $\overline{g} \in \argmin_{\overline{g} \in C(\overline{u})} \langle \overline{g}, -\overline{v} \rangle$ with $\lVert g - \overline{g} \rVert < \varepsilon$.
\label{lem:lin_opt_pert}
\end{mylemma}

\begin{proof}~\\
Let $( \overline{u}, \overline{v}) \in \H \times \H$. We can describe the set $\argmin_{\overline{g} \in C(\overline{u})} \langle \overline{g}, -\overline{v} \rangle$ using the vertices of $C(\overline{u})$. The set $C(\overline{u})$ is a convex polyhedron and the objective function $\langle \overline{g}, -\overline{v} \rangle$ is linear. A minimum of $\min_{\overline{g} \in C(\overline{u})} \langle \overline{g}, -\overline{v} \rangle$ is attained at a vertex of $C(\overline{u})$ and therefore it exists at least one $i \in \lbrace 1, \dots, m \rbrace$ such that $\langle \nabla f_i(\overline{u}), - \overline{v} \rangle = \min_{\overline{g} \in C(\overline{u})}\langle \overline{g}, -\overline{v} \rangle$. The same can be done for any $(u,v) \in \H \times \H$. Define the sets of optimal and nonoptimal vertices
\begin{align}
\begin{split}
    \overline{\mathcal{A}} &\coloneqq \lbrace i \in \lbrace 1,\dots, m \rbrace\, : \, \langle \nabla f_i(\overline{u}), -\overline{v} \rangle = \min_{\overline{g} \in C(\overline{u})}\langle \overline{g}, -\overline{v} \rangle \rbrace, \text{ and}\\
    \overline{\mathcal{I}} &\coloneqq \lbrace 1,\dots, m \rbrace \setminus \overline{\mathcal{A}}.\\
    \mathcal{A}(u,v) &\coloneqq \lbrace i \in \lbrace 1,\dots, m \rbrace\, : \, \langle \nabla f_i(u), -v \rangle = \min_{g \in C(u)}\langle g, -v \rangle \rbrace, \text{ and}\\
    \mathcal{I}(u,v) &\coloneqq \lbrace 1,\dots, m \rbrace \setminus \mathcal{A}(u,v).
\end{split}
\end{align}
There exists $M \in \R$ such that for all $i \in \overline{\mathcal{A}}$ and $j \in \overline{\mathcal{I}}$ it holds that
\begin{align}
    \langle \nabla f_i(\overline{u}), -\overline{v} \rangle < M < \langle \nabla f_j(\overline{u}), -\overline{v} \rangle.
\end{align}
Then by the continuity of $(u,v) \mapsto \langle \nabla f_i(u), v \rangle$ we can choose $\delta > 0$ such that for all $i \in \overline{\mathcal{A}}$ and $j \in \overline{\mathcal{I}}$
\begin{align}
    \langle \nabla f_i(u), -v \rangle < M < \langle \nabla f_j(u), -v \rangle,
\end{align}
for all $(u,v) \in \H \times \H$ with $\lVert (u,v) - (\overline{u}, \overline{v}) \rVert_{\H \times \H} < \delta$. For these $(u,v)$ it holds that $\mathcal{A}(u,v) \subset \overline{\mathcal{A}}$. Now, the rest follows from the continuity of the function $\nabla f_i(\cdot)$. Let $g \in \argmin_{g \in C(u)} \langle g, -v \rangle$. We write $g = \sum_{i \in \mathcal{A}(u,v)} \lambda_i \nabla f_i(u)$ as a convex combination of the optimal vertices in $C(u)$. Since $\mathcal{A}(u,v) \subset \overline{\mathcal{A}}$, it follows that $\overline{g} = \sum_{i \in \mathcal{A}(u,v)} \lambda_i \nabla f_i(\overline{u})$ is a solution to $\min_{\overline{g} \in C(\overline{u})} \langle \overline{g}, -\overline{v} \rangle$. Since all $\nabla f_i(\cdot)$ are continuous, we can choose $\delta > 0$ such that $\lVert g - \overline{g} \rVert < \varepsilon$.
\end{proof}

We are now in the position to prove Proposition \ref{prop:set_valued_G}.
\begin{proof}[Proof of Proposition \ref{prop:set_valued_G}:]~\\
\label{proof:prop_set_valued_G}

    (i) Fix $(t,u,v) \in [t_0, +\infty) \times \H \times \H$. The set $C(u) \coloneqq \co\left( \lbrace \nabla f_i(u) : i = 1,\dots, m\rbrace \right)$ is convex and compact as a convex hull of a finite set. Then $\argmin_{g \in C(u)} \langle g, -v \rangle$ is also convex and compact and the statement follows since sums and Cartesian products of convex and compact sets are convex and compact.\\

    (ii) We show that $G$ is u.s.c. in the $\varepsilon$ sense using Lemma \ref{lem:lin_opt_pert}. Then, we use Proposition \ref{prop:usc_varepsilon_sense} together with (i) to conclude $G$ is u.s.c. as well. This is technical but we include the proof for the sake of completeness.

    Using Lemma \ref{lem:lin_opt_pert} we can show that for all $\varepsilon > 0$ there exists $\delta > 0$ satisfying
    \begin{align*}
        G(B_{\delta}((\overline{t}, \overline{u}, \overline{v}))) \subset G(\overline{t}, \overline{u}, \overline{v}) + B_{\varepsilon}((0,0)),
    \end{align*}
    where $B_{\delta}((\overline{t}, \overline{u}, \overline{v}))) \subset [t_0, + \infty) \times \H \times \H$ and $B_{\varepsilon}((0,0)) \subset \H \times \H$ are open balls with radius $\delta$ and $\varepsilon$, respectively. To this end, we show that for all $(t, u, v) \in \H \times \H$ with $\lVert (\overline{t}, \overline{u}, \overline{v}) - (t, u, v) \rVert_{\R \times \H \times \H} < \delta$ and for all $(x,y) \in G(t,u,v)$ there exists an element $(\overline{x}, \overline{y}) \in G(\overline{t}, \overline{u}, \overline{v})$ with $\lVert (\overline{x}, \overline{y}) - (x, y) \rVert_{\H \times \H} < \varepsilon$. \\
    
    For $(t,u,v) \in \R \times \H \times \H$, $(x,y) \in G(t, u, v)$ is equivalent to
    \begin{align}
        \begin{split}
        x &= v, \\
        y &= - \frac{\alpha}{t} v -g , \text{ with}\\
        g &\in \argmin_{g\in C(u)} \langle g, -v \rangle.
        \end{split}
    \end{align}
    From Lemma \ref{lem:lin_opt_pert}, we know that there exists $\delta_1 > 0$ such that from $\lVert (t, u, v) - (\overline{t}, \overline{u}, \overline{v}) \rVert_{\H\times\H} < \delta_1$ it follows that there exists $\overline{g} \in \argmin_{\overline{g}\in C(\overline{u})} \langle \overline{g}, -\overline{v} \rangle$ such that
    \begin{align}
    \label{eq:eps_drittel_1}
        \lVert g - \overline{g} \rVert < \frac{\varepsilon}{3}.
    \end{align}
    Fix $\overline{g} \in \argmin_{\overline{g}\in C(\overline{u})} \langle \overline{g}, -\overline{v} \rangle$ satisfying \eqref{eq:eps_drittel_1}. Further, there exists $\delta_2 > 0$ such that from $\lvert t - \overline{t} \rvert < \delta_2$ it follows that
    \begin{align}
    \label{eq:eps_drittel_2}
        \left\lvert \frac{\alpha}{t} - \frac{\alpha}{\overline{t}} \right\rvert \lVert v \rVert < \frac{\varepsilon}{3}.
    \end{align}
    Let $\delta = \min\left\lbrace \delta_1, \delta_2, \frac{\varepsilon}{3(1+\alpha/t_0)}\right\rbrace$. It holds that
    \begin{align*}
        (\overline{x}, \overline{y}) \coloneqq \left(\overline{v}, - \frac{\alpha}{\overline{t}} \overline{v} - \overline{g} \right) \in G(\overline{t}, \overline{u}, \overline{v}).
    \end{align*}
    Then it follows that 
    \begin{align*}
        & \lVert (x,y) - (\overline{x}, \overline{y}) \rVert_{\H\times\H} \le \lVert v - \overline{v} \rVert + \left\lVert - \frac{\alpha}{t} v - g + \frac{\alpha}{\overline{t}} \overline{v} + \overline{g} \right\rVert \\
        \le & \left(1 + \frac{\alpha}{\overline{t}} \right) \lVert v - \overline{v} \rVert + \left\lvert \frac{\alpha}{\overline{t}} - \frac{\alpha}{t} \right\rvert \lVert v \rVert + \lVert g - \overline{g} \rVert < \varepsilon,
    \end{align*}
    which completes the proof. \\

    (iii) If $\dim(\H) < + \infty$ the proof follows from (ii). On the other hand, from $\phi$ being locally compact, we follow that $v \mapsto v$ is locally compact which is equivalent to $\H$ being finite-dimensional. \\

    (iv)  Before we start with the proof, we recall that the norm $\lVert (\cdot, \cdot) \rVert_{\H \times \H}$ and the norm $\lVert \cdot \rVert$ fulfill the following inequality. For all $x,y \in \H$, it holds that
    \begin{align*}
        \left\lVert (x,y) \right\rVert_{\H \times \H} \le \lVert x \rVert + \lVert y \rVert \le \sqrt{2} \left\lVert (x,y) \right\rVert_{\H \times \H}.
    \end{align*}
    Let $(t,u,v) \in [t_0, + \infty) \times \H \times \H$ and $\xi \in G(t,u,v)$. Then $\xi = \left( v, - \frac{\alpha}{t} v - g \right)$, with $g \in \argmin_{g \in C(u)} \langle g, -v \rangle$ and we follow
    \begin{align*}
        &\left\lVert \xi \right\rVert_{\H \times \H} \le \lVert v \rVert + \left\lVert  - \frac{\alpha}{t} v - g \right\rVert\\
        \le &\left(1 + \frac{\alpha}{t}\right) \lVert v \rVert + \max_{\theta \in \Delta^m} \left\lVert \sum_{i=1}^m \theta_i \nabla f_i(u) \right\rVert  \\
        \le & \left(1+ \frac{\alpha}{t_0}\right) \lVert v \rVert + \max_{\theta \in \Delta^m} \left\lVert \sum_{i=1}^m \theta_i \left(\nabla f_i(u) - \nabla f_i(0)\right) \right\rVert + \max_{\theta \in \Delta^m} \left\lVert \sum_{i=1}^m \theta_i \nabla f_i(0) \right\rVert \\
        \le & \left(1 + \frac{\alpha}{t_0}\right) \lVert v \rVert + L \lVert u \rVert + \max_{i=1,\dots,m} \left\lVert \nabla f_i(0) \right\rVert \\
        \le & c\left(1 + \left\lVert (u,v) \right\rVert_{\H \times \H} \right),
    \end{align*}
where we choose $c = \sqrt{2}\max\left\lbrace \left(1 + \frac{\alpha}{t_0}\right), L, \max_{i=1,\dots,m} \left\lVert \nabla f_i(0) \right\rVert \right\rbrace$.
\end{proof}

\end{document}